\documentclass[a4paper,11pt]{article}
\usepackage{amsmath,amscd,amssymb}
\usepackage[frenchb]{babel}
\usepackage{graphicx,psfrag}

\def \C {{\mathbb C}}
\def \l {{\lambda}}
\def \te {{\Theta_{\p_2,5}^{reg}}}
\def \U {{\mathfrak{U}}}
\def \M {{\mathfrak{M}}}
\def \A {{\mathfrak{A}}}
\def \ts {{\Theta_{sing}}}

\def \t {{\Theta}}
\def \cC {{\mathfrak{C}}}

\def \P {{\cal P}}

\def \p {{\mathbb P}}
\def \pu {{\mathbb P}_1}

\def \G {{\mathbb G}}
\def \fll {\longrightarrow}
\def \fl {\rightarrow}

\def \K {{{\cal K}}}
\def \L {{{\cal L}}}
\def \vs {\vskip}
\def \hilb {{\rm Hilb}_{\p_4,(8,5)}^{lisses, nd, nt}}
 \def \GA {{\Gamma}}
\def \hs {\hskip}

\def \I {{\cal I}}

\def \oo {{\cal O}}
\def \dm {{\textit{Preuve --- }}}

\def \a {{\alpha}}
\def \b {{\beta}}
\def \d {{\gamma}}

\def \vp {\varphi}
\def \opu {{{\cal O}_{\pu}}}

\def \pic {{\rm Pic}}
\def \hom {{\underline{\rm {\bf Hom}}}}
\def \tor {{\underline{\rm {\bf Tor}}}}

\def \noi {\noindent}

\def \mapdown#1{\Big\downarrow \rlap{$\vcenter{\hbox{$\scriptstyle#1$}}$}}

\def \Z {{\mathbb{Z}}}

\def \th #1#2{\noi\textsc{#1 #2.} ---}
\def \thg #1#2#3{\noi{\bf{#1}}{\bf\uppercase{{\footnotesize{#2}}} #3.} ---}
\def \thf #1#2{\noi\textsc{#1 #2.} ---}

\begin{document}

~
\vs -1.5 cm

\centerline{\large{\uppercase{\bf{Quelques remarques sur }}}}
\centerline{\large{\uppercase{\bf{les courbes de genre 5}}}}

\vs 0.4 cm

\centerline{\large{Jean \textsc{d'Almeida}, Laurent \textsc{Gruson} et
    Nicolas \textsc{Perrin}}}



\vs 0.8 cm

\centerline{\Large{\bf Introduction}}

\vs 0.8 cm

Dans ce texte, nous nous int{\'e}ressons {\`a} la g{\'e}om{\'e}trie des
courbes de genre 5 et plus pr{\'e}cisement aux courbes de genre 5 munies
d'une involution sans point fixe. Nous donnons deux caract{\'e}risa\-tions
g{\'e}om{\'e}triques de l'existence d'une telle involution.
%
La premi{\`e}re de ces
caract{\'e}ri\-sa\-tions, conjectur{\'e}e dans \cite{ACGH}, est le
point de d{\'e}part de notre {\'e}tude.


\vs 0.4 cm

Consid{\'e}rons $C$ une courbe de genre 5 non hyperelliptique et non
trigonale plong{\'e}e canoniquement dans $\p_4$. Les quadriques qui la
contiennent forment alors un r{\'e}seau (c'est-{\`a}-dire
que l'on a $h^0\I_C(2)=3$). Dans le plan de ce r{\'e}seau, les
quadriques singuli{\`e}res forment une courbe $\GA$ de degr{\'e}
5. La g{\'e}om{\'e}trie de $\GA$ et celle de $C$ sont li{\'e}es (cf. \cite{ACGH}
pp 270-274, nous rappelerons quelques r{\'e}sultats au premier
paragraphe).
En particulier, lorsque $C$ est un rev{\^e}tement
double non ramifi{\'e} d'une courbe de genre 3, la quintique $\GA$ se
d{\'e}compose en la r{\'e}union d'une conique et d'une cubique. Les auteurs de
\cite{ACGH} posent la question de la r{\'e}ciproque :
si la courbe $\GA$ est la r{\'e}union d'une conique et d'une cubique, la
courbe $C$ est-elle rev{\^e}tement double non ramifi{\'e} d'une courbe de
genre 3~? Nous r{\'e}pondons par l'affirmative {\`a} cette question au premier
paragraphe en montrant le th{\'e}or{\`e}me suivant,

\vs 0.4 cm

\thg{T}{h{\'e}or{\`e}me}{0.1} \textit{La courbe $C$ est rev{\^e}tement double
  d'une courbe de genre 3 si et seulement si $\GA$ est r{\'e}union d'une
  conique et d'une cubique,}

\vs 0.4 cm

\noi
compl{\'e}tant ainsi l'{\'e}tude de \cite{ACGH}. Nous rappelons
{\'e}galement, dans ce paragraphe, un isomorphisme birationnel entre
$\mathfrak{M}_5$ (l'espace des modules des courbes de genre 5) et
l'espace des modules des th{\'e}ta-caract{\'e}ristiques support{\'e}es
par une quintique plane.

\vs 0.6 cm

Dans le second paragraphe nous revenons au point de d{\'e}part de cette
{\'e}tude commune qui est un expos{\'e} de Jean d'Almeida au s{\'e}minaire
de g{\'e}om{\'e}trie alg{\'e}brique de Jussieu dans lequel il cherchait {\`a}
donner une r{\'e}ponse g{\'e}om{\'e}trique au probl{\`e}me pr{\'e}c{\'e}dent.
Cette {\'e}tude repose notamment sur la description des courbes
lisses de degr{\'e} 8 et de genre 5 de $\p_3$ qui ont une infinit{\'e} de
quadris{\'e}cantes. Nous reprenons en d{\'e}tail cette construction 
g{\'e}om{\'e}trique. Nous d{\'e}crivons une correspondance birationnelle entre
les courbes de genre 5 rev{\^e}tement double non ramifi{\'e} d'une courbe
de genre 3 et les quintiques planes r{\'e}union d'une conique et d'une
cubique modulo l'action du groupe d'Heisenberg $H$ associ{\'e} {\`a} la
situation (th{\'e}or{\`e}me 2.7). 
Cette construction permet notamment de donner une nouvelle
caract{\'e}risation g{\'e}om{\'e}trique des courbes de genre 5 non
hyperelliptiques rev{\^e}tement double non ramifi{\'e} d'une courbe de genre
3. Notons $\M_5^{i,i}$ (resp. $\M_5^{i,d}$) les courbes de genre 5
munies d'une involution sans point fixe dont la vari{\'e}t{\'e} de Prym
associ{\'e}e est ind{\'e}composable (resp. d{\'e}compos{\'e}e).

\vs 0.4 cm

\thg{T}{h{\'e}or{\`e}me}{0.2} \textit{Soit $C$ une courbe lisse de genre
  5.}

(\i) \textit{{\rm Premi{\`e}re caract{\'e}risarion :} On a $C\in\M_5^{i,i}$ si
  et seulement s'il existe un plongement ${\cal M}$  de degr{\'e} 8 de $C$
  dans $\p_3$ pour lequel la courbe $C$ a une infinit{\'e} de
  quadris{\'e}cantes.}

\textit{{\rm Seconde caract{\'e}risarion :} On a $C\in\M_5^{i,i}$ si
  et seulement s'il existe un plongement ${\cal M}'$  de degr{\'e} 7 de $C$
  dans $\p_3$ et une droite $L$ rencontrant $C$ en un point pour lesquels la courbe $C\cup L$ a une infinit{\'e} de quadris{\'e}cantes.}

\vs 0.1 cm

\textit{Dans cette situation, la courbe $Y$ des quadris{\'e}cantes est la m{\^e}me
  quelque soit le plongement, elle est lisse de genre 2 et telle que
  $J(Y)={\rm Prym}(C)$. Notons
  $\mathfrak{J}_0$ les diviseurs ${\cal M}$ de degr{\'e} 8 du premier
  type, $\mathfrak{J}_1$ les diviseurs ${\cal M}'$ de degr{\'e} 7 du
  second type et $\mathfrak{J}$ la r{\'e}union de ces ensembles.  Il y a
  un morphisme $\mathfrak{J}\to J(Y)$ qui est un fibr{\'e} principal
  homog{\`e}ne de groupe $H$ (le groupe d'Heisenberg de $Y$).}

\vs 0.1 cm

(\i\i) \textit{On a $C\in\M_5^{i,d}$ si et seulement s'il existe un
  morphisme de $C$ dans $\p_3$ tel que son image $\overline{C}$ est
  li{\'e}e {\`a} une droite $L$ par une intersection de deux c{\^o}nes cubiques,
  de degr{\'e} 8 et a deux points doubles aux sommets des c{\^o}nes.}


\vs 0.1 cm

\textit{Dans cette situation, notons $\mathfrak{J}$ les diviseurs de
  degr{\'e} 8 de $C$ qui d{\'e}finissent de tels morphismes. On a un morphisme
  $\mathfrak{J}\to E_1\times E_2$ qui est un fibr{\'e} principal homog{\`e}ne
  de groupe $H$ (les $E_i$ sont les courbes elliptiques d{\'e}finissant
  les c{\^o}nes cubiques et $H$ le groupe d'Heisenberg associ{\'e} {\`a}
  $E_1\times E_2$).}






\vs 0.2 cm

\noi
{\bf Remerciements} : Le troisi{\`e}me auteur remercie Atanas Iliev
pour les nombreuses r{\'e}f{\'e}rences qu'il lui a communiqu{\'e}es. 

\section{La construction de \cite{ACGH}}

Rappelons les principaux r{\'e}sultats d{\'e}crits dans 
\cite{ACGH} pp 270-274 : 
%
%
soit $C$ une courbe de genre 5 non hyperelliptique plong{\'e}e
canoniquement dans $\p_4$ (on identifiera $C$ {\`a} son image).
Les quadriques contenant $C$ forment un
r{\'e}seau, c'est-{\`a}-dire ${\rm dim}_{\C}(H^0\I_C(2))=h^0\I_C(2)=3$.

\subsection{Le cas trigonal}

La courbe $C$ est trigonale si et seulement si toutes les quadriques
du r{\'e}seau sont singuli{\`e}res. On ne peut dans ce cas d{\'e}finir la
quintique $\GA$ de l'introduction.

\vs 0.2 cm

Dans toute la suite nous supposerons que $C$ n'est pas trigonale. On
d{\'e}finit alors la courbe $\GA$ comme lieu des quadriques
singuli{\`e}res du r{\'e}seau. C'est une quintique plane.


\subsection{L'{\'e}tude g{\'e}n{\'e}rale}

Pour une courbe $C$ g{\'e}n{\'e}rale, on v{\'e}rifie que le r{\'e}seau de
quadriques est {\'e}galement g{\'e}n{\'e}ral et que la courbe $\GA$ est non
singuli{\`e}re. 
%
%
Les singularit{\'e}s de $\GA$ correspondent 
aux quadriques de
rang 3  du r{\'e}seau et $\GA$ n'a que des point doubles ordinaires
pour singularit{\'e}s (en particulier $\GA$ n'a pas de composante
multiple).

%
%
%
%
%
%

\subsubsection{Les droites contenues dans $\GA$}

Les droites contenues dans $\GA$ sont en bijection avec les pinceaux
bi-elliptiques de $C$ (c'est-{\`a}-dire les morphismes de degr{\'e} 2 de $C$
vers une courbe elliptique). Il y en a donc au plus 5.

%


\subsubsection{Les coniques de $\GA$}


On passe enfin au cas o{\`u} la courbe $\GA$ se d{\'e}compose en une conique
et une cubique. Supposons tout d'abord que $C$ est une rev{\^e}tement
double non ramifi{\'e} d'une courbe $X$ de genre 3. On montre alors que
$\GA$ est r{\'e}union d'une cubique $C_1$ et d'une conique $C_0$. 

%
%
%
%
%
%
%

\vs 0.4 cm

Pour donner une caract{\'e}risation des courbes de genre 5 munies d'une
involution sans point fixe, 
nous montrons la r{\'e}ciproque : 

%

\vs 0.2 cm

\thg{T}{h{\'e}or{\`e}me}{1.1} \textit{Si $\GA$ est r{\'e}union d'une conique
  $C_0$ et d'une cubique $C_1$, alors la courbe $C$ est rev{\^e}tement
  double non ramifi{\'e} d'une courbe $X$ de genre 3.}

\vs 0.4 cm

Donnons une autre pr{\'e}sentation classique du
probl{\`e}me par 
les r{\'e}seaux de quadriques de $\p_4$. Consid{\'e}rons 
l'espace projectif $\p({\rm
  Hom}(\C^3,S^2\C^5)^\vee)$ sur lequel agissent les groupes $PGL_3$ et
$PGL_5$. Dans cet espace projectif consid{\'e}rons l'ouvert $U$ 
des $\vp\in\p({\rm Hom}(\C^3,S^2\C^5)^\vee)$ tels que 
\begin{itemize}
\item[(\i)] $\vp:\C^3\to S^2\C^5$ est injective (pour qu'on ait un vrai
  r{\'e}seau de quadriques),
\item[(\i\i)] ${\rm det}(\vp)$ est non nul (vu comme d{\'e}terminant d'une matrice
  sym{\'e}trique de taille $5\times 5$ {\`a} coefficients dans
  ${\C^3}^\vee$),
\item[(\i\i\i)] aucune quadrique du r{\'e}seau n'est d{\'e}g{\'e}n{\'e}r{\'e}e
  en deux plans,
\item[(\i v)] l'intersection dans $\p_4$ des quadriques du r{\'e}seau
  forme une courbe lisse de genre 5. 
\end{itemize}

\vs 0.2 cm

\noi
On s'int{\'e}resse maintenant aux actions de $PGL_3$ et $PGL_5$ sur
$U$.

\vs 0.2 cm

\thg{F}{ait}{1.2} \textit{L'action de $PGL_3$ est libre et tous les
  points de $U$ sont stables pour cette action.} 

\vs 0.2 cm

\dm Ceci vient simplement de la condition (\i). Le quotient est
alors un ouvert de $\G(3,S^2\C^5)$ la grassmannienne des sous-espaces
vectoriels de dimension 3 de $S^2\C^5$.\hfill$\Box$

\vs 0.2 cm

\thg{L}{emme}{1.3} \textit{Les points de $U$ sont stables pour
  l'action de $PGL_5$.}

\vs 0.2 cm

\dm En effet, C.T.C Wall a {\'e}tudi{\'e} la stabilit{\'e} des
r{\'e}seaux de quadriques dans $\p_n$ (th{\'e}or{\`e}me 0.1 de
\cite{Wall}). En particulier, il montre que si un point $\vp$ est non
stable, alors la courbe de $\p_2$ d{\'e}finie par le d{\'e}terminant ${\rm
  det}(\vp)$ a une composante multiple (corollaire 2,
\cite{Wall}). Ceci ne peut jamais arriver dans notre situation : on
sait que la courbe $C$ de genre 5 intersection des quadriques du
r{\'e}seau est lisse et dans ce cas la courbe $\GA$ (cf. construction de
\cite{ACGH}) qui a pour {\'e}quation ${\rm det}(\vp)$ n'a pas de
composante multiple.\hfill$\Box$ 

\vs 0.2 cm

Les deux quotients g{\'e}om{\'e}triques $\U_3=U/PGL_3$ et $\U_5=U/PGL_5$
sont munis respectivement d'une action de
$PGL_5$ et $PGL_3$. Par un r{\'e}sultat de C. T. C. Wall \cite{wall2},
comme tous les points de $U$ sont stables pour les deux actions, ils
sont tous stables pour l'action de $PGL_3\times PGL_5$ et les deux
quotients $\U_3/PGL_5$ et $\U_5/PGL_3$ sont isomorphes au quotient
$U/(PGL_3\times PGL_5)$.

\vs 0.1 cm

Nous donnons deux interpr{\'e}tations de ces deux quotients. Notons
$\Theta_{\p_2,5}^{reg}$ l'ouvert de la vari{\'e}t{\'e} des
th{\'e}ta-caract{\'e}ristiques planes (cf. \cite{Sorger} pour son existence en
toute g{\'e}n{\'e}ralit{\'e}) support{\'e}ees par une courbe de degr{\'e} 5 
r{\'e}guli{\`e}res (c'est-{\`a}-dire $h^0\theta=0$).

\vs 0.2 cm

\thg{P}{roposition}{1.4} \textit{La vari{\'e}t{\'e} $\U_5$ est isomorphe {\`a}
  l'ouvert $\t_0$ de $\te$ correspondant aux conditions (\i\i\i) et
  (\i v).}

\vs 0.2 cm

\dm On a un morphisme naturel de $\Psi:U\to\te$ donn{\'e} de la mani{\`e}re
suivante : si $\vp\in U$, alors on d{\'e}finit un fibr{\'e} vectoriel sur
$\p_2=\p(\C^3)$ par le conoyau $\theta$ de la fl{\`e}che 
$$(\C^5)^\vee\!\otimes\oo_{\p_2}(-2)\stackrel{\vp}{\fl}\C^5\otimes
\oo_{\p_2}(-1)\ \ \ \ \ \ \ \ \ \ \ (*)$$ 
qui est {\'e}videment une th{\'e}ta-caract{\'e}ristique 
r{\'e}guli{\`e}re support{\'e}e par une courbe de degr{\'e} 5. Le morphisme est bien
d{\'e}fini sur tout $U$ et la fl{\`e}che $(*)$ 
est toujours injective car
le d{\'e}terminant ${\rm det}(\vp)$ n'est jamais nul. Ce morphisme est
{\'e}videment invariant sous $PGL_5$.

Ce morphisme est surjectif sur l'ouvert $\t_0$. En
effet, si $\theta\in\te$, alors la suite spectrale de Beilinson
(cf. \cite{OSS}) nous donne la suite spectrale suivante :
$$E_1^{i,j}=H^j\theta(1-i)\otimes\Omega_{\p_2}^i(i)\Longrightarrow
\theta(1)\ {\rm en}\ {\rm degr\acute{e}}\ 0.$$
Ici tous les termes sont nuls sauf $E_1^{0,0}$ et $E_1^{2, 1}$ ce qui
nous donne la r{\'e}solution :
$$0\fl H^1\theta(-1)\otimes\oo_{\p_2}(-1)\stackrel{\vp}{\fl}
H^0\theta(1)\otimes\oo_{\p_2}\fl\theta(1)\fl 0$$
o{\`u} $\vp$ est sym{\'e}trique car $\theta$ est une th{\'e}ta-caract{\'e}ristique,
donc on peut voir $\vp$ comme un {\'e}l{\'e}ment de $U$ (car
$\theta\in\t_0$) en fixant des bases duales de
$H^1\theta(-1)$ et $H^0\theta(1)$.

Il reste {\`a} montrer que les fibres sont les orbites pour l'action de
$PGL_5$ sur $U$. Or les {\'e}lements $\vp\in U$ de la fibre au-dessus de
$\theta\in\te$ sont les matrices sym{\'e}triques $5\times 5$ tel que le
quotient du morphisme de $\oo_{\p_2}$-module $(*)$ est $\theta$. On
peut alors identifier $\C^5$ {\`a} $H^0\theta(1)$ et $\vp$ est uniquement
donn{\'e}e {\`a} action de $PGL_5$ pr{\`e}s par la surjection
$H^0\theta(1)\otimes\oo_{\p_2}\fl\theta(1)$.\hfill$\Box$

\vs 0.2 cm

Notons $\hilb$ le sch{\'e}ma de Hilbert des courbes lisses de
degr{\'e} 8 et de genre $5$ de $\p_4=\p(\C^5)$ non
d{\'e}g{\'e}n{\'e}r{\'e}es (c'est-{\`a}-dire non contenues dans un
hyperplan) et non trigonales.

\vs 0.2 cm

\thg{P}{roposition}{1.5} \textit{On a un isomorphisme
  $\U_3\simeq\hilb$.} 

\vs 0.2 cm

\dm
On a un morphisme $\Phi:U\to\hilb$ invariant sous $PGL_3$ donn{\'e} par
l'intersection des quadriques du r{\'e}seau. Ce morphisme est bien {\`a}
valeur dans $\hilb$ gr{\^a}ce aux hypoth{\`e}ses faites sur $U$ (en
particulier comme toutes les quadriques ne sont pas
singuli{\`e}res --- i.e. ${\rm det}(\vp)\neq0$ --- la
courbe n'est pas trigonale). Il est surjectif car toute courbe lisse
$C$ de $\p_4$ de degr{\'e} 8 et de genre 5, non d{\'e}g{\'e}n{\'e}r{\'e}e est plong{\'e}e
canoniquement et $C$ est n{\'e}c{\'e}ssairement non hyperelliptique. On a
alors $h^0\I_C(2)=3$ ce qui donne un r{\'e}seau de quadriques $\vp\in U$
et la courbe est l'intersection des quadriques du r{\'e}seau.

Cette construction montre que les fibres au-dessus de $C\in\hilb$ sont
exactement donn{\'e}es par les orbites sous $PGL(H^0K_C)$ c'est-{\`a}-dire
sous $PGL_5$.\hfill$\Box$

\vs 0.2 cm

Notons $\M^0_5$ l'ouvert de l'espace des modules des courbes lisses de
genre 5 form{\'e} des courbes non hyperelliptiques et non trigonales. 

\vs 0.2 cm

\thg{C}{orollaire}{1.6} \textit{On a un isomorphisme
  $\M_5^0\simeq\t_0/PGL_3$.}

\vs 0.2 cm

\dm
Comme on sait que $\U_3/PGL_5\simeq\U_5/PGL_3$, il suffit de constater
que l'on a $\hilb/PGL_5\simeq\M_5^0$. Cette correspondance
birationnelle entre $\M_5$ et $\te/PGL_3$ est d{\'e}j{\`a} d{\'e}crite dans
\cite{ACGH}, ch VI, appendix C, th{\'e}or{\`e}me page 301. Nous montrons ici
simplement qu'elle est d{\'e}finie sur tout $\M_5^0$ c'est-{\`a}-dire d{\`e}s que
la courbe est non hyperelliptique et non trigonale. Remarquons par
ailleurs que cette description est utilis{\'e}e dans \cite{Kats} pour
montrer que la vari{\'e}t{\'e} $\M_5^i$ des courbes de genre 5 munies
d'une involution sans point fixe est rationnelle.\hfill$\Box$

\vs 0.2 cm

Nous revenons maintenant au th{\'e}or{\`e}me 1.1.

\vs 0.2 cm

\textit{Preuve du th{\'e}or{\`e}me 1.1.} --- Il s'agit de montrer que si $C$
est une courbe de genre 5 telle que la courbe $\GA$ (d{\'e}finie au
paragraphe 1.1) se d{\'e}compose en une conique $C_0$ et une cubique
$C_1$, alors la courbe $C$ est un rev{\^e}tement double non ramifi{\'e} d'une
courbe $X$ lisse de genre 3. Nous montrons gr{\^a}ce {\`a} la description
pr{\'e}c{\'e}dente le lemme suivant qui est le nouvel ingr{\'e}dient permettant la
d{\'e}monstration du th{\'e}or{\`e}me 1.1 :

\vs 0.2 cm

\thg{L}{emme}{1.7} \textit{Il existe un {\'e}l{\'e}ment $\vp\in\p({\rm
  Hom}(\C^3,S^2\C^5)^\vee)$ qui s'envoie sur $C\in\M_5^0$ par $\Phi$
et tel que dans des bases de $\C^5$ la matrice de $\vp$ peut s'{\'e}crire
sous la forme diagonale par blocs
$$\left(\begin{array}{cc}
 A & 0 \\
 0 & B
\end{array}\right)$$
o{\`u} $A$ est une matrice sym{\'e}trique de taille $2\times 2$ {\`a}
coefficients dans $(\C^3)^\vee$ (dont le d{\'e}terminant donne une
{\'e}quation de $C_0$) et $B$ est une matrice sym{\'e}trique de
taille $3\times 3$ {\`a} coefficients dans $(\C^3)^\vee$ (dont le
d{\'e}terminant donne une {\'e}quation de $C_1$).}

\vs 0.2 cm

\dm
Il s'agit donc de montrer que toutes les th{\'e}ta-caract{\'e}ristiques
$\theta\in\te$ ayant pour support $C_0\cup C_1$ peuvent {\^e}tre d{\'e}finie {\`a}
partir d'une matrice de cette forme. 

Soit $\theta$ une telle th{\'e}ta-caract{\'e}ristique et 
supposons que la conique $C_0$ est lisse. Nous savons (cf. 1.2) 
que les quadriques de rang 3 du r{\'e}seau
correspondent exactement aux points singuliers du support de $\theta$
(qui a seulement des points doubles ordinaires). La
fibre du faisceau $\theta$ est donc de dimension 0 en g{\'e}n{\'e}ral sur
$\p_2$, de dimension 1 pour un point g{\'e}n{\'e}ral de $C_0\cup C_1$ et de
dimension 2 aux six points d'intersection de $C_0$ et
$C_1$ (et {\'e}ventuellement en un point singulier de $C_1$). Notons $Z$
le sch{\'e}ma de cette intersection et consid{\'e}rons la restriction de la
r{\'e}solution~de~$\theta$
$$0\fl(\C^5)^\vee\otimes\oo_{\p_2}(-2)\stackrel{\vp}{\fl}
\C^5\otimes\oo_{\p_2}(-1)\fl\theta\fl 0$$
{\`a} la conique $C_0$. On a alors :
$$0\fl \tor_1^{\oo_{\p_2}}(\theta,\oo_{C_0})\fl
(\C^5)^\vee\otimes\oo_{C_0}(-2)\stackrel{\vp}{\fl}
\C^5\otimes\oo_{C_0}(-1)\fl\theta\vert_{C_0}\fl 0.$$
Nous identifions $C_0$ {\`a} $\pu$, il existe alors deux entiers $a$
et $b$ tels que la suite exacte pr{\'e}cedente s'{\'e}crive :
$$0\fl \opu(-b)\fl (\C^5)^\vee\otimes\opu(-4)\stackrel{\vp}{\fl}
\C^5\otimes\opu(-2)\fl\opu(a)\oplus\oo_Z\fl 0.$$ 
Cependant le morphisme central est sym{\'e}trique par hypoth{\`e}se donc si on
applique le foncteur $\hom_{\opu}(\ \bullet\ ,\opu(-6))$, on a la
m{\^e}me suite exacte. Ceci impose en particulier la relation $a=b-6$. Par
ailleurs un calcul de classes de Chern impose que $a+b=4$. Ainsi on a
$a=-1$ et $b=5$. 

Nous identifions maintenant le faisceau $Q$ image de la fl{\`e}che
$\vp$. C'est un faisceau localement libre de rang 4. Il s'{\'e}crit donc
$Q\simeq\oplus_{i=1}^4\opu(-a_i)$ avec $a_i$ des entiers compris entre 2
et 4. Par ailleurs, la sym{\'e}trie de la fl{\`e}che $\vp$ impose que $Q$
v{\'e}rifie la suite exacte suivante : 
$$0\fl Q\fl Q^\vee(-6)\fl\oo_Z\fl 0.$$
Ceci impose la condition $\sum_ia_i=15$ et on a 
$Q\simeq\opu(-3)\oplus\opu(-4)^3$. Il existe donc 
un sous-espace vectoriel de dimension 2 isomorphe {\`a} $H^1\opu(-3)$
contenu dans $(\C^5)^\vee$ et de fa\c con duale un quotient de rang 2
de $\C^5$ isomorphe {\`a} $H^0\opu(1)$.

Consid{\'e}rons la fl{\`e}che $Q\stackrel{M}{\fl}Q^\vee(-6)$ qui dans des
d{\'e}compositions de $Q$ et $Q^\vee(-6)$ est de la forme 
$$\left(\begin{array}{cc}
 \a & \b \\
 \ \!\!\!^t\b & \d
\end{array}\right)$$
o{\`u} $\a\in\C^*$ (nous verrons en fin de preuve que le cas $\a=0$ est
impossible), $\b$ est une matrice $1\times 3$ {\`a}
coefficients dans $H^0\opu(1)$ et $\d$ est une matrice $3\times 3$ {\`a}
coefficients dans $H^0\opu(2)$. En jouant sur les d{\'e}compositions de
$Q$ c'est-{\`a}-dire en faisant agir le groupe 
$${\rm Aut}(Q)=\left\{\left(\begin{array}{cc}
\lambda & u \\
0 &\mu
\end{array}\right)\ {\rm o\grave{u}}\ \lambda\in\C^*,\ \mu\in
GL_3(\C)\ {\rm et }\
u\ {\rm est}\ 
{\rm \grave{a}}\ {\rm coefficients}\ {\rm dans}\ H^0\opu(1)\right\},$$
on peut trouver une d{\'e}composition dans laquelle $M$ s'{\'e}crit sous la
forme
$$\left(\begin{array}{cc}
 1 & 0 \\
 0 & \d'
\end{array}\right)\ \ \ {\rm o\grave{u}}\ \d'\ {\rm est}\ {\rm une}\
{\rm matrice}\ 3\times3\  {\rm \grave{a}}\ {\rm coefficients}\ {\rm
  dans}\ H^0\opu(2).$$
La fl{\`e}che $\vp$ 
peut {\^e}tre {\'e}crite de la mani{\`e}re suivante :
$$\opu(-4)^5\fl Q\stackrel{M}{\fl}Q^\vee(-6)\fl\opu(-2)^5.$$
On fixe la d{\'e}composition pr{\'e}cedente de $Q$, il existe une
d{\'e}composition de $\C^5$ sous la forme $H^0\opu(1)\oplus\C^3$ (et donc
de $(\C^5)^\vee$ sous la forme $H^1\opu(-3)\oplus\C^3$) telle que la
fl{\`e}che de $\opu(-4)^5$ dans $Q$ soit de la forme suivante :
$$\left(\begin{array}{cc}
 m & 0 \\
 0 & I_3
\end{array}\right)$$
o{\`u} $m$ est l'application canonique de $\opu(-4)\otimes
H^1\opu(-3)\to\opu(-3)$ et $I_3$ est la matrice identit{\'e} de taille
$3\times 3$. La fl{\`e}che de $Q^\vee(-6)$ dans $\opu(-2)^5$ est la
transpos{\'e}e de cette fl{\`e}che. La compos{\'e}e est alors dans ces
d{\'e}compositions de la forme :
$$\left(\begin{array}{cc}
 m\ \!^tm & 0 \\
 0 & \d'
\end{array}\right)$$
o{\`u} $m$ et $\d'$ sont les matrices d{\'e}finies ci-dessus. C'est une matrice
{\`a} coefficients dans 
$$H^0\opu(2)\simeq H^0\oo_{\p_2}(1)=(\C^3)^\vee$$ 
et nous avons bien la forme diagonale par blocs recherch{\'e}e.

Si le morphisme $\a:\opu(-3)\to\opu(-3)$ est nul, alors la conique est
singuli{\`e}re. En effet le m{\^e}me raisonnement donne alors un matrice de la
forme 
$$\left(\begin{array}{cc}
 0 & \!\!^tm\ \!\b \\
 \!^t\!\b\ \!m & \d'
\end{array}\right).$$
Mais $\b:\opu(-4)^3\to\opu(-3)$ peut se mettre sous la forme
$\b=(0,\b_1,\b_2)$ avec $\b_1$ et $\b_2$ dans $H^0\opu(1)$. Si on
{\'e}crit $\b'=(\b_1,\b_2)$ on a une matrice de la forme (ici $s\in
H^0\opu(2)$) :
$$\left(\begin{array}{ccc}
 0 & 0 &\!\!^tm\ \!\b \\
 0 & s & * \\
 \!^t\!\b\ \!m & * &* 
\end{array}\right)$$
dont le d{\'e}terminant a $\det(\ \!\!^tm\ \!\b)$ comme facteur double. C'est
impossible.

Dans le cas o{\`u} la conique est singuli{\`e}re, on raisonne de la m{\^e}me
mani{\`e}re en se restreignant aux deux droites. Dans ce cas, la
matrice $A$ peut en plus se diagonaliser.\hfill$\Box$

\vs 0.2 cm

Nous terminons maintenant la preuve du th{\'e}or{\`e}me 1.1. La courbe $C$ de
genre 5 est l'intersection d'un r{\'e}seau de quadriques diagonalisables
par blocs :
$$\left(\begin{array}{cc}
 A & 0 \\
 0 & B
\end{array}\right).$$
Nous munissons $C$ d'une involution : consid{\'e}rons la
sym{\'e}trie de $\p_4$ donn{\'e}e par la matrice
$$\left(\begin{array}{cc}
 -I_2 & 0 \\
 0 & I_3
\end{array}\right),$$
elle laisse le r{\'e}seau et la courbe $C$ invariants. Il reste {\`a}
v{\'e}rifier que cette involution est sans point fixe. Sur $C$ les
points fixes 
sont les points d'intersection de la courbe avec la droite ou le plan
d{\'e}finis par la d{\'e}composition de $\C^5$. Les points fixes sont donc
donn{\'e}s par les points base du r{\'e}seau $A$ de quadriques de $\pu$ ou du
r{\'e}seau $B$ de quadriques planes. Pour le r{\'e}seau $A$, il n'y a pas de
tels points car la conique d{\'e}finie par le d{\'e}terminant serait
alors une droite double. 
Pour le r{\'e}seau $B$, la formule d'H{\"u}rwitz nous dit que le
nombre de points fixes est multiple de 4. Or un r{\'e}seau de coniques a au
plus 3 points base, l'involution est donc sans point fixe.\hfill$\Box$

\vs 0.2 cm

Nous pouvons alors compl{\'e}ter le tableau de \cite{ACGH} page 274 qui
r{\'e}sume les liens entre la g{\'e}om{\'e}trie de $C$ et celle de $\GA$ ($\ts$ est
le lieu singulier du diviseur th{\'e}ta de la jacobienne de $C$) :
\begin{itemize}
\item $C$ est trigonale $\Leftrightarrow$ $\GA=\p_2$ $\Leftrightarrow$
  $\t_{sing}$ est la r{\'e}union de deux copies de $C$ conjugu{\'e}es sous
  l'involution et qui se rencontrent en deux points,
\item $C$ n'est pas trigonale $\Leftrightarrow$ $\GA$ est une
  quintique plane $\Leftrightarrow$
  $\t_{sing}/-1=\GA$,
\item $C$ est bi-elliptique $\Leftrightarrow$ $\GA$ contient une droite
  $\Leftrightarrow$ $\t_{sing}$ contient une composante elliptique qui
  est un rev{\^e}tement double de la droite,
\item $C$ est rev{\^e}tement double d'une courbe de genre 3
  $\Leftrightarrow$ $\GA$ contient une conique $\Leftrightarrow$
  $\t_{sing}=\Sigma_1\cup\Sigma_2$ les deux composantes se coupant en
  six points avec $p_a(\Sigma_1)=2$ et $p_a(\Sigma_2)=4$,
\item $C$ n'est pas trigonale et n'a pas de pinceau semi-canonique
  $\Leftrightarrow$ $\GA$ est une quintique plane lisse
  $\Leftrightarrow$ $\t_{sing}$ est irr{\'e}ductible et lisse de genre 11.
\end{itemize}

%

\vs 0.2 cm

Notons $\t_0^{3,2}$ le ferm{\'e} de $\t_0$ form{\'e} des
th{\'e}ta-caract{\'e}ristiques dont le support est la r{\'e}union d'une conique et
d'une cubique et notons $\M_5^{i}$ le localement ferm{\'e} de l'espace
des modules des courbes de genre 5 form{\'e} des courbes non
hyperelliptiques et non trigonales qui sont munies d'une involution
sans point fixe. Le th{\'e}or{\`e}me 1.1 et le corollaire 1.6 permettent de
montrer le

\vs 0.2 cm

\thg{C}{orollaire}{1.8} \textit{On a un isomorphisme entre
  $\M_5^{i}$ et $\t_0^{3,2}$.}

\vs 0.2 cm

\thf{Remarques}{1.9} (\i) Ce corollaire permet de red{\'e}montrer le
r{\'e}sultat de P.I. Katsylo \cite{Kats} sur la rationnalit{\'e} de $\M_5^i$. En
effet, $\t_0^{3,2}$ est le produit de la vari{\'e}t{\'e} des coniques planes
qui est {\'e}videment rationnelle par la vari{\'e}t{\'e} des cubiques planes munies d'une
th{\'e}ta-caract{\'e}ristique paire qui est elle aussi rationnelle
(cf. \cite{DolKan} corollaire 5.7.2).

(\i\i) Nous allons dans la suite {\'e}tudier une autre
description de $\M_5^i$. Nous remarquons ici que si $C$ est une courbe
lisse de genre 5 non hyperelliptique munie d'une involution sans point
fixe $i$, alors $C$ est non trigonale. 

En effet, supposons que $C$ soit trigonale. Soit $D$ un $g^1_3$ sur
$C$, si $D$ est invariant par $i$ (i.e. $i^*D=D$), alors l'involution
$i$ induit une involution sur $\pu$ par le morphisme d{\'e}duit de
$D$. Cette derni{\`e}re a n{\'e}cessairement au moins un point
fixe. La fibre du morphisme au-dessus d'un point fixe {\'e}tant stable
par $i$ et form{\'e}e de trois points, elle a au moins un point fixe par
$i$ ce qui est absurde. Si $D$ n'est pas fix{\'e} par l'involution, soit
$D'=i^*D$. On peut alors regarder le morphisme 
$$C\to\p(H^0D)\times\p(H^0D')=\pu\times\pu\subset\p(H^0D\otimes H^0D')=\p_3.$$
L'image de $C$ est une courbe de bidegr{\'e} $(3,3)$ dans la quadrique,
c'est-{\`a}-dire une courbe de genre arithm{\'e}tique $4$, c'est impossible
car $C$ est de genre 5.

\section{Une tentative d'interpr{\'e}tation g{\'e}om{\'e}trique}

Dans ce paragraphe nous d{\'e}crivons d'une nouvelle mani{\`e}re l'espace de
module des courbes de genre 5 rev{\^e}tement double d'une courbe de genre
3 gr{\^a}ce aux quintiques planes r{\'e}union d'une conique et d'une cubique.

\vs 0.2 cm

Remarquons tout d'abord que si $C$ est une courbe de genre 5
rev{\^e}tement double d'une courbe de genre 3, elle est alors munie d'une
involution $i$ sans point fixe. 
La jacobienne de $C$ est, elle aussi, munie d'une involution toujours
not{\'e}e $i$ et on note ${\rm Prym}(C,i)={\rm Im}(i-1)$ la vari{\'e}t{\'e} de
Prym associ{\'e}e (cf. par exemple \cite{ACGH} ch. VI appendix C ou
\cite{mumprym}). C'est alors une surface ab{\'e}lienne principalement
polaris{\'e}e. On a alors un morphisme 
$${\rm Prym}:\M_5^i\to\A_2$$
o{\`u} $\A_2$ est l'espace de modules (grossier) des surfaces ab{\'e}liennes
principalement polaris{\'e}es. L'image d'une courbe $C$ est en g{\'e}n{\'e}ral une
surface ab{\'e}lienne ind{\'e}composable, c'est alors la jacobienne d'une
courbe $Y$ de genre 2. Nous notons $\M_5^{i,i}$ l'ouvert de $\M_5^{i}$
correspondant {\`a} ce cas. Dans le cas contraire,  la vari{\'e}t{\'e} de Prym est
le produit de deux courbes elliptiques et nous notons $\M_5^{i,d}$ le
ferm{\'e} correspondant (le premier cas correspond {\`a} une conique
irr{\'e}ductible contenue dans $\GA$ alors que dans le second cas la
conique est d{\'e}compos{\'e}e en deux droites).


\vs 0.2 cm

Nous commen\c cons par d{\'e}crire la fibre de ce morphisme. Elle a {\'e}t{\'e}
{\'e}tudi{\'e}e par de nombreux auteurs (voir par exemple \cite{Don},
\cite{Dona}, \cite{Donsmi} et \cite{Verra}). En particulier A. Verra
\cite{Verra} d{\'e}crit compl{\`e}tement cette fibre (en incluant les
rev{\^e}tement admissibles introduits par \cite{Beauvilleprym}).
Nous incluons une preuve de ses r{\'e}sultats dans le cas des rev{\^e}tements
de courbes lisses et decrivons plus en d{\'e}tail le cas d'un produit de
courbes elliptiques.

Nous aurons besoin de quelques rappels sur les courbes de genre 2,
leur jacobienne et la surface de Kummer associ{\'e}e. Pour plus de d{\'e}tails,
voir \cite{Hud} ou \cite{mum}.

\subsection{Notations et rappels}

Soit $A\in\A_2$ une surface ab{\'e}lienne principalement polaris{\'e}e et soit
$\t_A$ son diviseur th{\'e}ta. Ce diviseur d{\'e}finit un morphisme 
$$A\to\p(H^0(2\t_A))\simeq\p_3.$$
L'image $\K$ de ce morphisme est de quotient $A/{\rm Aut}(\t_A)$.

Si $A$ est ind{\'e}composable cette image est la surface de Kummer
$\K\simeq J(Y)/-1$ associ{\'e}e {\`a} la courbe $Y$ de genre 2 telle que
$J(Y)=A$. Si $A$ est d{\'e}composable isomorphe au produit de courbes
elliptiques $E_1\times E_2$, alors le morphisme est de degr{\'e} 4 et
l'image est une quadrique lisse~de~$\p_3$.

Le groupe $H$ des {\'e}l{\'e}ments d'ordre 2 de $A$ est fini d'ordre 16
isomorphe {\`a} $(\Z/2\Z)^4$. Ce groupe est le quotient du groupe
d'Heisenberg ${\cal H}$ {\`a} 32 {\'e}l{\'e}ments par son centre (cf. \cite{mum}
ou \cite{mumtata3}). Par abus de notation nous l'appelerons encore
groupe d'Heisenberg associ{\'e} {\`a} $A$. Il agit lin{\'e}airement sur $\p_3$ et
laisse stable la surface $\K$. 

Notons $Z$ son image dans $\K$. Dans le
premier cas $Z$ est l'ensemble des points doubles de $\K$ et $H$ est
le groupe des automorphismes de $\K$ (ils laissent n{\'e}cessairement fixe
$Z$ voir \cite{Hud} ou \cite{mum} p. 353). Dans le second cas, on a
une involution $i_k$ donn{\'e}e par le diviseur $2\t_{E_k}$ sur chacune
des courbes $E_k$. Le lieu de ramification sur $\K$ du
morphisme $A\to\K$ est la r{\'e}union de 4 droites de chaque
famille. L'ensemble $Z$ est le groupe de 16 points obtenus comme
intersection de ces droites. Le groupe $H$ est le sous-groupe du
groupe des automorphismes de $\K$ laissant stable $Z$. 

Remarquons que l'on a une involution $i_1\times i_2$ sur $A$. Notons
$\widetilde{\K}$ le quotient de $A$ par cette involution. On a alors
une suite de morphismes de degr{\'e} 2 $A\to\widetilde{\K}\to\K$.
Le premier est ramifi{\'e} au dessus de l'image r{\'e}ciproque de $Z$ dans
$\widetilde{\K}$ alors que le second est ramifi{\'e} au dessus des quatres
droites de chaque famille. 

Notons $U(A)$ l'ouvert de $\p_3^\vee$ --- invariant sous $H$
--- des plans de $\p_3$ qui ne sont pas tangents {\`a} la surface $\K$ et
ne passant pas par un point de $Z$. Notons $\K^\vee\subset\p_3^\vee$
la surface de duale de $\K$ et $Z^\vee$ la r{\'e}union des 16 plans
duaux des 16 points de $Z$. L'ouvert $U(A)$ est $\p_3^\vee$ priv{\'e}
de $\K^\vee\cup Z^\vee$.

\subsection{Lien entre $C$ et $A$}

\thg{P}{roposition}{2.1} \textit{On a un isomorphisme
$$U(A)/H\simeq{\rm Prym}^{-1}(A).$$}
%
\vs -0.5 cm

\dm Nous commen\c cons par associer {\`a} tout point de $x\in U(A)$ une
courbe lisse de genre 5 rev{\^e}tement double d'une courbe de genre
3. 
Le point $x$ est dans $\p_3^\vee$, il d{\'e}finit donc un plan
$x^\vee$ de $\p_3$ qui 
rencontre donc $\K$ en une courbe lisse $\widetilde{X}$ (car $x\in U(A)$).

Dans le premier cas, la courbe $X=\widetilde{X}$ est une quartique
lisse donc de genre 3. Consid{\'e}rons la courbe $C\subset A$ image
r{\'e}ciproque de $\widetilde{X}$ par le morphisme $A\to\K$. La courbe $C$
est un rev{\^e}tement double de $X$. Il est non ramifi{\'e} car $X$ {\'e}vite les
points doubles de $\K$. 

Dans le second cas, la courbe $\widetilde{X}$ est une conique
lisse. Consid{\'e}rons encore la courbe $C\subset A$ image r{\'e}ciproque de
$\widetilde{X}$ par le morphisme $A\to\K$. On peut par ailleurs
relever $\widetilde{X}$ en ${X}$ dans $\widetilde{\K}$. Ce rel{\`e}vement
est ramifi{\'e} en 8 points. La courbe ${X}$ est donc lisse de genre
3. Comme $x\in U(A)$ le morphisme induit de $A\to\widetilde{\K}$ donne
un rev{\^e}tement non ramifi{\'e} $C\to X$.

La classe d'isomorphisme de $C$ ne d{\'e}pend pas du repr{\'e}sentant $x$ de
$\overline{x}$ : si on change de repr{\'e}sentant, la courbe $C$ est
translat{\'e}e par un {\'e}l{\'e}ment $\a$ d'ordre 2 dans $A$ et la courbe
$\widetilde{X}$ est obtenue par la transformation lin{\'e}aire
correspondante dans $\p_3$. Nous verrons en fin de d{\'e}monstation
pourquoi $C$ est non hyperelliptique.

\vs 0.2 cm

R{\'e}ciproquement soit $C\in{\rm Prym}^{-1}(A)$ et $X$ le quotient de $C$
par l'involution.
Fixons $\xi$ un diviseur de degr{\'e} 1 sur $C$ tel que $2\xi=2i(\xi)$ (il
suffit de prendre $\xi$ tel que $2\xi=x_0+i(x_0)$ avec
$x_0\in C$, ce qui existe toujours car une vari{\'e}t{\'e} ab{\'e}lienne est
divisible). Consid{\'e}rons le morphisme
$$C\to A\ \ \ \ \ \ \ \ $$
\vs -0.9 cm
$$\ \ \ \ \ \ \ \ \ \ \ \ \ \ \ x\mapsto(1-i)(x-\xi)\ \ \ \ \ \ \ \ $$
obtenu en composant l'injection de $C$ dans $J(C)$ donn{\'e}e par
$x\mapsto x-\xi$ avec la fl{\`e}che $1-i$ de $J(C)$ dans $A$.

\vs 0.2 cm

\thg{L}{emme}{2.2} \textit{Le morphisme $C\to A$ est une immersion.}

\vs 0.2 cm

\dm
Nous commen\c cons par montrer l'injectivit{\'e} de la fl{\`e}che. Soient
$x$ et $y$ dans $C$ tels que leurs images soient les m{\^e}mes dans
$A$. On a alors $x-i(x)=y-i(y)$. Posons $D=x+i(y)$, c'est un
diviseur de degr{\'e} 2 sur $C$ et on a $D=i(x)+y$. Si $x\neq y$, le
diviseur $D$ a au moins 2 sections (car les paires $\{x,i(y)\}$ et
$\{i(x),y\}$ ne peuvent {\^e}tre {\'e}gales). Dans ce cas la courbe $C$ a un
diviseur de degr{\'e} 2 avec au moins 2 sections donc $C$ est
hyperelliptique. C'est absurde d'o{\`u} l'injectivit{\'e} de la fl{\`e}che.


Il reste {\`a} prouver que l'application tangente est {\'e}galement
injective. Nous avons vu que le morphisme $C\fl A$ se factorise par
la fl{\`e}che $1-i:J(C)\to A$. Consid{\'e}rons l'application
$$f:C\to\p(T_{0}J(C))$$
d{\'e}finie par : 
soit $x\in C$ et soit $v_x$ un
vecteur tangent {\`a} $C$ en $x$, alors $f(x)$ est le translat{\'e} en
$0$ de l'image de $v_x$ par la diff{\'e}rentielle de l'application
$C\to J(C)$ d{\'e}finie par $x\mapsto x-\xi$. 

L'application $f$ est exactement le plongement canonique de $C$ (voir
par exemple \cite{ACGH} exercices A-5 et A-6 p. 263). Le noyau en
chaque point de la diff{\'e}rentielle de 
$$J(C)\to A$$ 
$$\ \ \ \ \ \ \ \ \ \xi\ \ \mapsto\ (1-i)(\xi)$$
est donn{\'e} par l'espace tangent de la jacobienne $J(X)$ de
$X$. L'application tangente de $C\to A$ est injective si et seulement
si la courbe $C$ dans son plongement
canonique ne rencontre pas l'espace projectif correspondant {\`a}
$\p(T_0J(X))$. Ce dernier est le plan projectif associ{\'e} au
sous-espace propre pour la valeur propre 1 de l'involution $i$ sur
$C$.
La courbe $C$ ne rencontre
jamais ce plan car l'involution $i$ est sans point fixe. 
L'application tangente est donc injective.\hfill$\Box$

\vs 0.2 cm

La courbe $C$ est donc plong{\'e}e dans $A$ et est invariante par
l'involution $x\mapsto -x$ de $A$ qui correspond {\`a} l'involution $i$
sur $C$. Soit alors $\widetilde{X}$ l'image de $C$ dans
$\K\subset\p_3$. Il nous suffit de montrer que $\widetilde{X}$ est une
section hyperplane (de degr{\'e} 4 ou 2 selon les cas). Cependant comme le
degr{\'e} du morphisme est 2 ou 4 selon les cas, il suffit de montrer que  
%
$$C\cdot(1-i)^*(2\t_A)=8.$$

Fixons quelques notations, si $J$ est une vari{\'e}t{\'e} ab{\'e}lienne
principalement polaris{\'e}e, nous notons $\widehat{J}$ sa duale, $\t$ son
diviseur th{\'e}ta et $\l_\t:J\to\widehat{J}$ l'isomorphisme associ{\'e}. 
Dans notre situation, on a un morphisme
$J(X)\stackrel{\pi^*}{\to}J(C)$ et on peut d{\'e}finir 
$$\psi=\l_{\t_X}^{-1}\circ\widehat{\pi^*}\circ\l_{\t_A}:J(C)\to J(X).$$
On a {\'e}galement le morphisme d'inclusion $j:A\to J(C)$ et le
morphime $1-i:J(C)\to A$. On peut alors d{\'e}finir les morphismes 
$$\sigma=\pi^*\times j:J(X)\times A\to J(C)\ \ {\rm et}\ \
\tau=\psi\times(1-i):J(C)\to J(X)\times A$$
et on a (cf. \cite{mumprym}) $\sigma\circ\tau=2\cdot{\rm Id}_C$.

\vs 0.2 cm

\thg{F}{ait}{\!\!} \textit{On a \hs 0.5 cm$4\t_C=2(1-i)^*\t_A+2\psi^*\t_X$.}

\vs 0.2 cm

\dm
Dans \cite{mumprym}, D. Mumford montre que le diagramme suivant
commute :
$$\begin{array}{ccc}
J(X)\times A & \stackrel{\l_{2\t_X}\times\l_{2\t_A}}{\fll} &
\widehat{J(X)}\times \widehat{A} \\
\mapdown{\sigma} & & \mapdown{\widehat{\sigma}} \\
J(C) & \stackrel{\l_{\t_C}}{\fll} & \widehat{J(C)}
\end{array}.$$
La compos{\'e}e $J(C)\stackrel{\tau}{\fl}J(X)\times A
 \stackrel{\l_{2\t_X}\times\l_{2\t_A}}{\fll} \widehat{J(X)}\times
\widehat{A} \stackrel{\widehat{\tau}}{\fl}\widehat{J(C)}$ est donc
$J(C)\stackrel{\l_{4\t_C}}{\fll} \widehat{J(C)}$.\hfill$\Box$

\vs 0.2 cm

On peut 
calculer notre intersection. En effet, on sait que
$C\cdot\t_C=g(C)=5$. Par ailleurs, la restriction de $\psi$ {\`a} $C$ est le
morphisme $\pi$ donc $C\cdot\psi^*\t_X=2X\cdot\t_X=2g(X)=6$. Ainsi on a 
$$C\cdot(1-i)^*(2\t_A)=4C\cdot\t_C-C\cdot2\psi^*\t_X=8.$$

\vs 0.2 cm

Il reste {\`a} montrer que dans la premi{\`e}re construction la courbe $C$
n'est pas hyperelliptique. Si elle {\'e}tait hyperelliptique, son image
par $f$ serait alors une courbe rationnelle normale $C_4$ de degr{\'e} 4
de $\p_4$. L'involution descendrait alors sur $C_4$ 
(la courbe $C_4$ est l'image du plongement canonique et est donc
invariante par $i$). Mais alors la courbe rationnelle rencontrerait
l'espace propre associ{\'e} {\`a} la valeur propre 1. Ceci impose que
l'application tangente de $C$ vers $A$ serait non injective en un
point. La courbe $C\in A$ aurait alors un cusp ce qui est
absurde.\hfill$\Box$
 
\vs 0.2 cm

Au-dessus de 
$\A_2$, on a un fibr{\'e} projectif de dimension 3
donn{\'e} par $\p(H^0(2\t_A)^\vee)$ au-dessus de $A$. On note
$\mathfrak{U}(\A_2)$ l'ouvert donn{\'e} dans chaque fibre au-dessus de $A$
par $U(A)$. Au-dessus de $\A_2$, on a aussi une fibration
$\mathfrak{H}$ en groupes d'ordre 16 correspondant au-dessus de $A$ au
groupe d'Heisenberg $H$ des {\'e}l{\'e}ments d'ordre 2 de
$A$. 
On a une action fibre {\`a} fibre de $\mathfrak{H}$ sur
$\mathfrak{U}(\A_2)$.

\vs 0.2 cm

\thg{C}{orollaire}{2.3} \textit{La vari{\'e}t{\'e} $\M_5^i$ des courbes de
  genre 5 lisses non hyperelliptiques et munies d'une involution sans
  point fixe est isomorphe au quotient
  $\mathfrak{U}(\A_2)/\mathfrak{H}$.}


\subsection{Lien avec les courbes planes de degr{\'e} 5}

Dans la situation pr{\'e}c{\'e}dente et pour $C\in\M_5^{i,i}$, la
vari{\'e}t{\'e} ab{\'e}lienne $A={\rm Prym}(C)$ est la jacobienne d'une courbe $Y$
de genre 2. Choisissons $H'$ un des 16 plans tangents {\`a} $\K^\vee$
selon une conique. Le plan $H'$ est tangent {\`a} $\K$ le
long de la conique $C_2$ image du plongement de Veronese de
$\p(H^0(K_Y))$. Ce plan est donc est isomorphe {\`a} $\p(S^2H^0(K_Y))$. La
conique $C_2$ contient 6 points doubles de la surface $\K^\vee$ (que
nous notons $x_k$ pour $k\in[1,6]$).

%
%

\vs 0.2 cm

\thg{P}{roposition}{2.4} \textit{La donn{\'e}e d'un point $x\in\p_3^\vee$
  est {\'e}quivalente {\`a} la donn{\'e}e d'une cubique de $H'$ passant par les six
  points $x_k$ pour $k\in[1,6]$.}

\vs 0.2 cm

\dm
{\`A} un point $x\in\p^\vee_3$ nous pouvons associer la
surface cubique $W$ polaire de $\K^\vee$. Cette surface cubique
contient en particulier les points doubles de $\K^\vee$. Son
intersection avec le plan $H'$ d{\'e}finit une cubique $C'$ qui passe par
les 6 points $x_k$ pour $k\in[1,6]$.

Nous d{\'e}crivons cette application. On prend des coordonn{\'e}es $(X,Y,Z,T)$
dans $\p_3^\vee$. Soit $F$ l'{\'e}quation de $\K^\vee$ et $X$ l'{\'e}quation
de $H'$. Nous {\'e}crivons alors
$$F=\lambda X^4+X^3F_1+X^2F_2+XF_3+F_4$$
o{\`u} les $F_i$ sont homog{\`e}nes de degr{\'e} $i$ en $(Y,Z,T)$. Par ailleurs
comme le plan $H'$ est tangent {\`a} $\K^\vee$ le long de la conique
$C_2$, on a $F_4=C_2^2$ (on a ici encore not{\'e} $C_2$ l'{\'e}quation de la
conique $C_2$). On {\'e}crit alors le point $x$ sous la forme $(a,b,c,d)$,
la surface $W$ a pour {\'e}quation
$$a\frac{\partial F}{\partial X}+b\frac{\partial F}{\partial Y}+
c\frac{\partial F}{\partial Z}+d\frac{\partial F}{\partial T}.$$
L'{\'e}quation de la courbe $C'$ dans le plan $H'$ est alors 
$$aF_3+b\frac{\partial F_4}{\partial Y}+c \frac{\partial F_4}{\partial
  Z}+d\frac{\partial F_4}{\partial T}.$$
L'application qui {\`a} $x$ associe $C'$ est donc lin{\'e}aire. Elle va d'un
espace projectif de dimension 3 dans un autre. Nous montrons qu'elle
est bijective : 
son image est de
dimension 3. 

Si le point $x$ est dans le plan $H'$, alors la courbe $C'$ a pour
{\'e}quation (rappelons que $F_4=C_2^2$) :
$$2C_2\bigg(b\frac{\partial C_2}{\partial Y}+c \frac{\partial C_2}{\partial
  Z}+d\frac{\partial C_2}{\partial T}\bigg).$$
La courbe $C'$ contient dans ce cas la
conique $C_2$ et on obtient un plan de cubiques de cette mani{\`e}re. Par
ailleurs lorsqu'on prend $x=(1,0,0,0)$, l'{\'e}quation de $C'$ est alors
simplement $F_3$. Cette cubique ne contient pas $C_2$ car sinon la
surface $\K^\vee$ serait singuli{\`e}re de long de toute la conique
$C_2$. L'image est donc de dimension 3.\hfill$\Box$

%

\vs 0.2 cm

\thf{Remarque}{2.5} On a ainsi 16 isomorphismes lin{\'e}aires --- $\phi_{H'}$
pour chaque plan $H'$ --- entre le $\p^\vee_3$ contenant $\K^\vee$ et
l'espace projectif de dimension 3, notons le $\p(N(Y))$, des cubiques du
plan $\p(S^2(H^0K_Y))$
passant par les six points $(x_k)_{k\in[1,6]}$. Chaque isomorphisme
permet de d{\'e}finir une action de $H$ 
sur $\p(N(Y))$ par $\phi_{H'}\sigma\phi_{H'}^{-1}$ pour $\sigma\in
H$. Toutes ces actions sont les m{\^e}mes car on a
$\phi_{\sigma(H')}=\phi_{H'}\circ\sigma^{-1}$ pour $\sigma\in H$. On a
donc une action bien d{\'e}finie de $H$ sur $\p(N(Y))$

%
%
%

\vs 0.2 cm

\thf{Remarque}{2.6} La vari{\'e}t{\'e} des courbes de genre 5 non
hyperelliptiques munies d'une involution sans point fixe et de vari{\'e}t{\'e}
de Prym $A=J(Y)$ est isomorphe (proposition 2.1) {\`a} l'ouvert
$U(A)/H$ de $\p_3^\vee/H$. Elle est donc isomorphe {\`a} un ouvert de
$\p(N(Y))/H$. Les orbites de cubiques de cet ouvert ne contiennent en
particulier pas la conique comme composante.

Par ailleurs, la donn{\'e}e de la conique canonique de $\p(S^2(H^0K_Y))$
et des six points $(x_k)_{k\in[1,6]}$ permet de retrouver la courbe
$Y$ comme rev{\^e}tement double de la conique ramifi{\'e} aux 6 points. 

La proposition 2.4 permet de r{\'e}interpr{\'e}ter l'ouvert $U(A)$ comme un
ouvert de l'espace $\p(N(Y))$ des cubiques du plan $\p(S^2(H^0K_Y))$
passant par les six points $x_k$ pour $k\in[1,6]$. Nous pouvons donc
r{\'e}interpr{\'e}ter le fibr{\'e} $\mathfrak{U}(\A_2)$ au-dessus de
$\A^i_2$ (les surfaces ab{\'e}liennes ind{\'e}compos{\'e}es) comme un ouvert
${\bf U}$ de la vari{\'e}t{\'e} des courbes planes de degr{\'e} 5
r{\'e}union d'une conique lisse et d'une cubique la coupant en 6
points distincts. Nous avons encore une action sur ${\bf U}$ du
fibr{\'e} en groupe d'Heisenberg $\mathfrak{H}$ au-dessus de $\A^d_2$.

\vs 0.2 cm

\thg{T}{h{\'e}or{\`e}me}{2.7} \textit{Les constructions pr{\'e}c{\'e}dentes d{\'e}crivent
  un isomorphisme entre la vari{\'e}t{\'e} $\M_5^i$ des courbes de genre 5 non
  hyperelliptiques munies d'une involution sans point fixe et la
  vari{\'e}t{\'e} ${\bf  U}/\mathfrak{H}$.}

\textit{Les orbites des courbes ne contiennent en particulier pas les
  courbes de degr{\'e} 5 contenant une conique double.} 

\vs 0.2 cm

\thf{Remarque}{2.8} Il existe une autre param{\'e}trisation (voir par
exemple \cite{L}) de la vari{\'e}t{\'e} $\M_5^i$ par les courbes planes de
degr{\'e} 5 r{\'e}union d'une conique $A$ et d'une cubique $B$ (sans
th{\'e}ta-caract{\'e}ristique). Cette param{\'e}trisation fait intervenir des
constructions classique de g{\'e}om{\'e}trie plane et est fortement li{\'e}e {\`a} la
construction de $\cite{ACGH}$ : les coniques $C_0$ et $A$ sont les
m{\^e}mes, la cubique $C_1$ est la Hessienne de la cubique $B$ et la
courbe $X$ est l'enveloppe d'un pinceau de coniques obtenu {\`a} partir de
$A$ et $B$.

\section{Des plongements particuliers}

Nous allons ici r{\'e}aliser l'isomorphisme de la proposition 2.1 de
mani{\`e}re g{\'e}om{\'e}trique, c'est-{\`a}-dire avec des courbes
plong{\'e}es. Nous montrerons ainsi que si $C\in\M_5^i$, alors elle
poss{\`e}de une famille de plongements particuliers qui caract{\'e}risent
l'existence de l'involution sans point fixe.


\subsection{Pr{\'e}liminaires}

\subsubsection{Le cas $\M_5^{i,i}$}

Soit $Y$ une courbe lisse de genre 2, on note $A$ sa jacobienne et
$\K$ la surface de Kummer associ{\'e}e. Fixons $\L$ un fibr{\'e}
inversible de degr{\'e} 6 sur $Y$. On consid{\`e}re le plongement de
$C$ donn{\'e} par $\L$
dans $\p_4=\p(H^0\L)$. On a toujours $h^0\I_Y(2)=4$
(cf. \cite{GLP}). On a donc un $\p_3$ 
de quadriques
contenant $Y$. Remarquons que comme il n'y a pas de diviseur $\xi$ de
degr{\'e} 3 sur $Y$ tel que $h^0\xi=3$, la courbe $Y$ n'a jamais de
tris{\'e}cante dans $\p_4$. 

Notons $\cC$ la sous-vari{\'e}t{\'e} de $\p_3$ des quadriques
singuli{\`e}res. C'est une vari{\'e}t{\'e} d{\'e}terminantielle dont la dimension
attendue est 2 et le degr{\'e} attendu 5. 


\vs 0.2 cm

\thg{P}{roposition}{3.1} \textit{La vari{\'e}t{\'e} $\cC$ est la r{\'e}union d'un
  plan et d'une surface de degr{\'e} 4 isomorphe {\`a} $\K^\vee$. Le plan est
  le plan dual du point correpondant {\`a} $\L\in J(Y)$.}

\vs 0.2 cm

\dm
Nous d{\'e}crivons toutes les quadriques $Q$ contenant la courbe
$Y$. 

\vs 0.1 cm

1. Si $Q$ est de rang inf{\'e}rieur ou {\'e}gal {\`a} 2, alors la quadrique est
   form{\'e}e de la r{\'e}union de deux hyperplans ou d'un hyperplan double. La
   courbe $Y$ est alors d{\'e}g{\'e}n{\'e}r{\'e}e ce qui est absurde.

\vs 0.2 cm

2. Si $Q$ est de rang 3, alors son noyau est une droite $L$ de
   $\p_4$. Cette droite rencontre $Y$ en au plus 4 points sinon la
   courbe $Y$ serait contenue dans un hyperplan (engendr{\'e} par la
   droite et 2 autres points de $Y$). Soit donc $x$ la longueur du
   sch{\'e}ma $Y\cap L$, on a $x=0,\ 2$ ou $4$ et consid{\'e}rons la
   projection de centre la droite $L$. Elle d{\'e}finit un morphisme de
   $Y$ vers une conique.

\vs 0.2 cm

\hs 0.1 cm 2.1. Si $x=4$, alors le morphisme est un isomorphisme ce
qui est absurde ($Y$ est de genre~2).

\vs 0.2 cm

\hs 0.1 cm 2.2. Si $x=2$, le mophisme est 
de degr{\'e} 2 et est 
donn{\'e} par la compos{\'e}e
du rev{\^e}tement double $Y\to\pu$ d{\'e}fini par $K_Y$ 
avec le plongement de Veronese
$\pu\to\p_2$. 
On a donc $\L=2K_Y+D$ o{\`u} $D$
est le diviseur correspondant aux deux points de $L\cap Y$. Le
diviseur $D$ est alors d{\'e}termin{\'e} et de degr{\'e} 2. Les deux points de
$L\cap Y$ sont donc uniquement d{\'e}termin{\'e}s sauf si $D=K_Y$. 

\vs 0.2 cm

\hs 0.6 cm 2.2.1. Ainsi, si $\L\neq 3K_Y$, alors $L$ est unique
(correspondant au diviseur effectif associ{\'e} {\`a} $\L-2K_Y$) et la
quadrique est {\'e}galement unique.

\vs 0.1 cm

\hs 0.6 cm 2.2.2. Si par contre $\L=3K_Y$, alors on a un $\pu$ de
choix pour la droite $L$ (toutes les bis{\'e}cantes qui coupent $Y$ en
deux points en involution), il y a donc le m{\^e}me $\pu$ de quadriques.

\vs 0.2 cm

\hs 0.3 cm 2.3. Si $x=0$, on a alors un morphisme de degr{\'e} 3 de $C$
vers $\pu$. Soit $D$ le diviseur associ{\'e},  on a ${\rm deg}D=3$ et
$h^0D=2$. Le morphisme vers la conique est donn{\'e} par 
$Y\to\p(H^0D)\stackrel{v}{\to}\p(S^2H^0D)$
o{\`u} $v$ est le plongement de Veronese. On a ici $\L=2D$. Si
$D_0$ est un tel diviseur, les autres  
diviseurs sont donn{\'e}s par $D+M$ avec $M$ une demi-p{\'e}riode
(i.e. $2M=0$). On a donc un  
nombre fini de tels diviseurs. La quadrique est alors
d{\'e}termin{\'e}e par $D$. La droite $L$ est le  
quotient de rang 2 suivant : $S^2(H^0D)\to H^0(2D)$ et la
conique l'image de $\p(H^0D)$ dans  
$\p(S^2H^0D)$. 

\vs 0.3 cm

3. Il reste le cas g{\'e}n{\'e}ral o{\`u} la quadrique est de rang 4, le noyau de
   $Q$ est alors un point $P$ et on projette par rapport {\`a} ce point
   (notons $\pi_P$ cette projection). On obtient alors un morphisme de
   $Y$ vers une quadrique de $\p_3$. Le morphisme $\pi_P\vert_Y$ est de degr{\'e}
   6 si $P\not\in Y$ et de degr{\'e} 5 sinon.

\vs 0.2 cm

\hs 0.3 cm 3.1. Supposons que $P\in Y$, la courbe $\pi_P(Y)$ est de
degr{\'e} divisant 5 dans $\p_3$. Elle ne peut  
{\^e}tre de degr{\'e} 1 sinon $Y$ serait d{\'e}g{\'e}n{\'e}r{\'e}e. Elle est donc
de degr{\'e} 5 (et de genre 2) et est contenue dans une unique
quadrique. Ceci nous donne une famille de dimension un~de~quadriques.

\vs 0.2 cm

\hs 0.3 cm 3.2. Reste le cas g{\'e}n{\'e}ral o{\`u} le point $P$ n'est pas dans
$Y$. L'image $\pi_P(Y)$ est alors une courbe  
de degr{\'e} divisant 6, c'est {\`a} dire de degr{\'e} 1, 2, 3 ou 6. De
plus $\pi_P(Y)$ ne peut {\^e}tre plane  
sinon $Y$ serait d{\'e}g{\'e}n{\'e}r{\'e}e. Donc $\pi_P(Y)$ est de degr{\'e} 3
ou 6 irr{\'e}ductible et n'est pas plane. 

\vs 0.2 cm

\hs 0.6 cm 3.2.1. Si ${\rm deg}(\pi_P(Y))=3$, alors la courbe
$\pi_P(Y)$ est une cubique gauche et le morphisme $Y\to\pi_P(Y)$ 
est de degr{\'e} 2 donn{\'e} par le diviseur $K_Y$ et la composition avec  
le plongement de Veronese $\p(H^0K_Y)\to\p(S^3H^0K_Y)$. On a
alors n{\'e}cessairement $\L=3K_Y$  
et le point $P$ est d{\'e}termin{\'e} par le conoyau de
$S^3H^0K_Y\to H^0(3K_Y)$. Dans ce cas $\pi_P(Y)$  
est contenue dans un plan de quadriques qui se rel{\`e}ve en un
plan dans $\p(H^0\I_Y(2)^\vee)$. 

\vs 0.1 cm

\hs 0.6 cm  3.2.2. Si ${\rm deg}(\pi_P(Y))=6$, 
la projection $Y\to\pi_P(Y)$ est birationnelle. La courbe $\pi_P(Y)$
est de degr{\'e} 6 contenue dans une quadrique lisse  
$\overline{Q}$ et son genre g{\'e}om{\'e}trique est 2. Dans
$\overline{Q}$, la classe de $\pi_P(Y)$ est $(a,b)$ avec $a+b=6$ et  
$(a-1)(b-1)\geq 2$. On a donc deux possibilit{\'e}s :
$(a,b)=(2,4)$ ou $(a,b)=(3,3)$. 

\vs 0.1 cm

Si $(a,b)=(2,4)$, la projection sur le premier facteur donne un
morphisme $\vp$ de degr{\'e} 2 de $Y$ vers $\pu$ (donc d{\'e}fini par
$K_Y$). 
L'image
r{\'e}ciproque d'un point du premier facteur de $\overline{Q}$ dans $C$
est une bis{\'e}cante en involution. La quadrique $Q$ contient donc la
surface $S$ recouverte par les bis{\'e}cantes en involution. Cette surface
est de degr{\'e} 3. Une quadrique contenant $Y$ contient $S$ si et
seulement si elle contient en plus un point $x_0$ de $S-C$. En effet,
dans ce cas on a $C\cup \{x_0 \}\subset S\cap Q$. Mais comme $Q$ est
de degr{\'e} 2 et $S$ de degr{\'e} 3, on a $Q\cap S$ est une courbe de degr{\'e} 6
ou $S\subset Q$. C'est donc que $S\subset Q$.  

L'ensemble des quadriques contenant $x_0$ forme un plan dans
$\p(H^0\I_Y(2))$. Par ailleurs toutes ces quadriques sont
singuli{\`e}res. En effet, si une de ces quadriques {\'e}tait lisse, la
surface $S$ serait alors intersection compl{\`e}te (sur une quadrique
lisse de $\p_4$, toute surface est intersection compl{\`e}te, car le
groupe de Picard de la quadrique est $\Z$). C'est absurde car son
degr{\'e} (qui est 3) devrait alors {\^e}tre un multiple de celui de $Q$ (qui
est 2). Ceci nous donne le plan.  Remarquons que si $\L=3K_Y$, la
surface $S$ est un c{\^o}ne et le plan pr{\'e}c{\'e}dent a {\'e}t{\'e} d{\'e}crit en 3.2.1. 

\vs 0.1 cm

Si $(a,b)=(3,3)$, la courbe $\pi_P(Y)$ a deux
points doubles. On a donc deux bis{\'e}cantes passant par $P$. Soient
$\xi_1$ et $\xi_2$ les diviseurs de degr{\'e}s 2 d{\'e}finis par ces
bis{\'e}cantes. Comme elles sont concourrantes, on a
$h^0(\L-\xi_1-\xi_2)\geq 2$. Or ${\rm deg}(\L-\xi_1-\xi_2)=2$ donc
$\L-\xi_1-\xi_2=K_Y$.  

Ainsi pour tout couple de points de $Y$, on plut{\^o}t pour tout diviseur
$\xi$ de degr{\'e} 2 sur $Y$, il existe un unique diviseur
$\xi'=\L-K_Y-\xi$ tel que les deux bis{\'e}cantes (si $\xi\neq\xi'$)
d{\'e}finies par des sections de $\xi$ et $\xi'$ se rencontrent en un
point $P$. Il y a ici deux cas particuliers a traiter {\`a} part : lorsque
$\xi$ (ou $\xi'$) est $K_Y$ et lorsque $\xi=\xi'$.  
 
En dehors de ces deux cas, la donn{\'e}e du couple $(\xi,\xi')$ d{\'e}finit
deux bis{\'e}cantes distinctes. Le point
d'intersection $P$ est donc bien d{\'e}fini. La projection de $Y$ par $P$
donne une courbe de degr{\'e} 6 ayant deux points doubles. Elle est donc
contenue dans une unique quadrique. 

Si $(\xi,\xi')=(K_Y,\L-2K_Y)$ ou sym{\'e}triquement si
$(\xi,\xi')=(\L-2K_Y,K_Y)$, on a un $\pu$ de bis{\'e}cantes correspondant
aux section de $K_Y$. Ces bis{\'e}cantes rencontrent toutes la droite
d{\'e}finie par $\L-2K_Y$ et par projection par rapport {\`a} cette droite on
retrouve les cas 2.2.1 et 2.2.2 o{\`u} la quadrique est de rang 3 et o{\`u}
$x=2$. Ainsi pour ces couples, on a alors un $\pu$ de quadriques si
$\L=3K_Y$ et une unique quadrique si $\L\neq K_Y$. 

Si enfin $\xi=\xi'$ (ie. si on a un point double de $\K^\vee$), les
deux bis{\'e}cantes sont alors confondues et les tangentes aux points
d'intersection de cette bis{\'e}cante avec $Y$ se rencontrent. On peut
alors choisir pour $P$ n'importe quel point de la bis{\'e}cante. La
projection aura alors un tacnode et sera donc toujours contenue dans
une unique quadrique. On a un $\pu$ de
quadriques.

\vs 0.1 cm

Ainsi, un couple de diviseurs $(\xi,\xi')$ d{\'e}finit un unique c{\^o}ne sauf
pour les points en involution sur la jacobienne : les $(\xi,\xi')$
pour lesquels on a un $\pu$ de quadriques. La surface obtenue est donc
isomorphe {\`a} la duale $\K^\vee$ de la surface de Kummer $\K$ de $Y$. 

\vs 0.3 cm

Montrons que le plan est bien le plan dual correspondant au
point $\L\in{\rm Pic}_6(Y)$. D{\'e}terminons l'intersection de ce plan et
de la surface $\K^\vee$. 

Si $\L=3K_Y$, 
l'intersection recherch{\'e}e est form{\'e}e des quadriques du plan
d{\'e}crit en 3.2.1 qui sont de rang 3. Elles forment une conique
double. Ce plan est donc le plan dual du point double de $\K$
qui correspond {\`a} $\L$. Ces quadriques sont exactement les quadriques
singuli{\`e}res dont le sommet rencontre $Y$.  

Si $\L\neq 3K_Y$, on a deux cas. Soit la quadrique a pour sommet une
droite et on est dans la situation 2.2.1. Soit la quadrique a pour
sommet un point et il doit appartenir {\`a} $Y$.
On voit que chaque point de la courbe $Y$ d{\'e}finit un point
de l'intersection du plan et de $\K^\vee$, les deux points du diviseur
$\L-2K_C$ se contractant en un m{\^e}me point. L'intersection du
plan et de $\K^\vee$ est donc une courbe de degr{\'e} 4 avec un point
double dont la d{\'e}singularisation est $Y$. C'est le plan dual du point
$\L$.\hfill$\Box$ 

\vs 0.2 cm 

\th{Remarque}{3.2} Soit $\L$ un diviseur de degr{\'e} 6 
sur $Y$ une courbe de genre 2. 
\textit{Se donner un point de $\p_3^\vee$
correspond exactement {\`a} se donner un fibr{\'e} vectoriel de rang 2
sur $Y$ de d{\'e}terminant~$\L$.}

En effet, se donner un point du $\p_3^\vee$ contenant $\K^\vee$
revient
{\`a} se donner une quadrique contenant la courbe $Y$ dans le plongement
$\p(H^0\L)$ (cf. proposition 3.1). Comme sur une quadrique de $\p_4$ il
y a un fibr{\'e} vectoriel de rang 2 tautologique, on a alors un fibr{\'e}
vectoriel de rang 2 et de d{\'e}terminant $\L$ sur $Y$. 

R{\'e}ciproquement si on a un fibr{\'e} vectoriel $E$ de rang 2 sur $Y$ et de
d{\'e}terminant $\L$, alors on a une forme quadratique sur $\Lambda^2H^0E$
qui induit une forme quadratique sur $H^0(\Lambda^2E)=H^0\L$. Cette
forme quadratique correspond {\`a} une quadrique contenant $Y$ et donc {\`a}
un point de $\p_3^\vee$.

Cette remarque permet de red{\'e}montrer un r{\'e}sultat de Narasimhan et
Ramanan d{\'e}crivant l'espace des modules des fibr{\'e}s vectoriels de
rang 2 et degr{\'e}s pair \cite{NR}

%
%
%
%
%
%

\subsubsection{Le cas $\M_5^{i,d}$}

Soit $A=E_1\times E_2$ une surface ab{\'e}lienne d{\'e}compos{\'e}e o{\`u}
les $E_i$ sont des courbes elliptiques avec un point marqu{\'e}. Soit
$Y=E_1\cup E_2$ la r{\'e}union transversale des deux courbes de genre
1 se rencontrant en leur point marqu{\'e}.
Notons $\K$ la surface de Kummer associ{\'e}e et fixons
$\L=\L_1\boxtimes\L_2$ un  diviseur de bidegr{\'e} $(3,3)$ sur $Y$. On
plonge alors $Y$ par $\L$ dans $\p(H^0\L)=\p_4$. On a toujours
$h^0\I_Y(2)=4$ et donc un $\p_3=\p(H^0\I_Y(2)^\vee)$ de quadriques
contenant $Y$.


Notons $\Pi_i$ le plan de $\p_4$ contenant $E_i$. 
%
%
Toutes les quadriques $Q$ contenant $Y$ contiennene les plans
$\Pi_i$. Elles sont donc singuli{\`e}res 
et leur sommet contient le point $P=\Pi_1\cap\Pi_2$.

Notons $\cC$ la sous-vari{\'e}t{\'e} 
des quadriques de rang inf{\'e}rieur {\`a} 3.

\vs 0.2 cm

\thg{P}{roposition}{3.3} \textit{La vari{\'e}t{\'e} $\cC$ est
  isomorphe {\`a} la surface $\K^\vee$.}

\vs 0.2 cm

\dm
On effectue la projection {\`a} partir du point $P$ et on arrive dans
un espace projectif de dimension 3. Les courbes $E_1$ et $E_2$ (ainsi
que les plans $\Pi_1$ et $\Pi_2$) ont pour images deux droites
disjointes.
On cherche toutes les quadriques de rang inf{\'e}rieur 
{\`a} 3 contenant ces deux droites.

Soit $Q$ une telle quadrique et supposons qu'elle est de rang
exactement 3. Soit alors $x$ son sommet (c'est un point). On projette
{\`a} partir de $x$, l'image de $Q$ doit {\^e}tre une conique lisse et
doit contenir l'image des deux droites. Ceci n'est jamais possible. La
quadrique $Q$ est donc r{\'e}union de deux plans chacun contenant l'une
des deux droites. 

Dans $\p_4$ les quadriques contenant $Y$ sont les r{\'e}unions
de deux hyperplans l'un contenant $\Pi_1$, l'autre contenant
$\Pi_2$. Ils forment une surface isomorphe {\`a} $\pu\times\pu$. On a un
morphisme de $E_1\times E_2$ vers cette surface : {\`a} un point de
$E_1$ on associe l'hyperplan contenant $\Pi_2$ et ce point et de
mani{\`e}re sym{\'e}trique pour $E_2$. On retrouve le recouvrement
$E_1\times E_2\to\K^\vee$.\hfill$\Box$

\subsection{La courbe $C$ de genre 5 rev{\^e}tement double d'une courbe
  $X$ de genre $3$.} 

Nous traitons les cas de $\M_5^{i,i}$ et $\M_5^{i,d}$ en m{\^e}me
temps. Soit $Y$ la courbe de genre 2 pr{\'e}c{\'e}dente (lisse ou r{\'e}union de
deux courbes elliptiques se coupant en un point) et $A$ la surface
ab{\'e}lienne correspondante. Soit $\L$ le diviseur de de degr{\'e} 6 qui
permet de plonger $Y$ dans $\p_4$, on choisit $Q$ une quadrique de
rang maximal contenant $Y$ (lisse dans le premier cas et de rang 4
dans le second). Elle correspond {\`a} un point $x$ du $\p_3^\vee$
contenant $\K^\vee$. Nous construisons ici de mani{\`e}re plong{\'e}e
l'isomorphisme d{\'e}crit {\`a} la proposition 2.1.

%

\vs 0.4 cm

Nous d{\'e}finissons la courbe $C$. Notons $A^\times$ l'ouvert des $\xi\in
A$ tels que $\xi\neq K_Y$ (si $Y$ est lisse) et $\xi\neq(P,y)$ ou
$(x,P)$ sinon. Alors $\xi\in A^\times$ d{\'e}finit une bis{\'e}cante $L_\xi$ de
$Y$. 
%
Consid{\'e}rons la courbe $C\subset A$ d{\'e}finie de la mani{\`e}re
suivante : 
$$C=\overline{\{\xi\in A^\times\ /\ L_\xi\subset Q\}}.$$
%
%
%

\vs 0.2 cm

\thg{L}{emme}{3.4} \textit{La courbe $C$ est invariante sous
  l'involution $x\mapsto -x$ de $A$.} 

\vs 0.2 cm

\dm
Dans le premier cas l'involution est donn{\'e}e par $\xi\mapsto\L-K_Y-\xi$.
Soit $\xi\neq K_Y$ et supposons que $\xi\in C$, alors $L_\xi\subset
Q$. Consid{\'e}rons maintenant $\xi'=\L-K_Y-\xi$, la bis{\'e}cante
$L_{\xi'}$ rencontre $L_\xi$. Mais alors $L_{\xi'}$ rencontre $Q$ en
trois points (les deux points de $L_{\xi'}\cap Y$ et le point de
$L_{\xi'}\cap L$). Elle est donc contenue dans $Q$ et on a $\xi'\in C$.

Dans le second cas l'involution est donn{\'e}e sur chacune des courbes
$E_i$ par $x\mapsto y$ o{\`u} $y$ est le troisi{\`e}me point
d'intersection de la droite $(xP)$ avec la courbe $E_i$ dans le plan
$\Pi_i$ ($P$ est le point d'intersection des courbes $E_i$). Si $L$
est une bis{\'e}cante {\`a} $Y$ contenue dans $Q$, la bis{\'e}cante $L'$ obtenue
par l'involution {\`a} partir de $L$ est contenue dans le plan
$\Pi=(P,L)$. Mais $P$ est dans le sommet de $Q$ et $L\subset Q$ donc
$\Pi\subset Q$ et donc $L'\subset Q$. La courbe $C$ est m{\^e}me
invariante par chacunes des involutions des courbes
$E_i$.\hfill$\Box$.

\vs 0.2 cm

\th{Remarque}{3.5} On a pas a priori donn{\'e} de condition pour que
$\xi\in A^\times$ appartienne~{\`a}~$C$. 

Dans le premier cas, si $\xi=K_Y$, on a un $\pu$ de sections et
donc autant de bis{\'e}cantes. Ces bis{\'e}cantes sont soit toutes en
m{\^e}me temps contenues dans $Q$ soit jamais contenues dans $Q$. 
Si $\L=3K_Y$, 
toutes les bis{\'e}cantes d{\'e}finies par $K_Y$ se coupent en un
m{\^e}me point. Ainsi le raisonnement pr{\'e}c{\'e}dent permet de
montrer que si l'une est dans $Q$ elles le sont toutes. Si par contre
$\L\neq 3K_Y$, alors les bis{\'e}cantes d{\'e}finies par $K_Y$
rencontrent toutes la m{\^e}me bis{\'e}cante $L_0$ d{\'e}finie par le
diviseur $\L-2K_Y$. Ainsi si une de ces bis{\'e}cante est dans $Q$,
alors $L_0$ est contenue dans $Q$ et par le raisonnement du lemme 3.4
elles sont toutes contenues dans $Q$.

Dans le second cas, $\xi=(x,P)$ (resp. $(P,y)$) est dans $C$ si et
seulement si le plan engendr{\'e} par la tangente en $P$ {\`a} la courbe
$E_2$ (resp. $E_1$) et le point $x$ (resp. $y$) est contenu dans $Q$.

\vs 0.4 cm

Soit $\xi\in A$ et consid{\'e}rons $L_\xi$ la bis{\'e}cante de $Y$
ainsi d{\'e}finie. On peut 
alors associer {\`a} $L_\xi$ le plan $H_\xi$ de
$\p(H^0\I_Y(2)^\vee)=\p_3^\vee$ des quadriques contenant $Y\cup
L_\xi$. C'est bien un plan car pour qu'une quadrique contenant $Y$
contienne $L_\xi$ il suffit qu'elle contienne un point de plus de
$L_\xi$. Le raisonnement de la remarque 3.5 permet de d{\'e}finir $H_\xi$
m{\^e}me pour $\xi\in A^\times$.

\vs 0.2 cm

\thg{L}{emme}{3.6} \textit{Le morphisme $A\to\p_3$ d{\'e}fini par
  $\xi\mapsto H_\xi$ est le morphisme de $A$ dans $\K$.}

\vs 0.2 cm

\dm
Le raisonnement du lemme 3.4 permet de montrer que ce morphisme est
invariant par l'involution sur $A$ (et m{\^e}me par celles des $E_i$ si
$Y$ n'est pas lisse). 

Ce morphisme est bien {\`a} valeur dans $\K$ car les plans de quadriques
ainsi d{\'e}finis sont tangents {\`a} la vari{\'e}t{\'e} $\K^\vee\subset\cC$ de
quadriques singuli{\`e}res. En effet, dans des bases bien choisie les
matrices des quadriques s'{\'e}crivent (dans chacun des cas) sous la forme 
$$\left(\begin{array}{ccccc}
0 & 0 & 0 & x & y \\
0 & 0 & a & a_{2,4} & a_{2,5} \\
0 & a & 0 & a_{3,4} & a_{3,5} \\
x & a_{2,4} & a_{3,4} & a_{4,4} & a_{4,5} \\
y & a_{2,5} & a_{3,5} & a_{4,5} & a_{5,5} 
\end{array}\right)
\ \ \ \ {\rm ou}\ \ \ \ 
\left(\begin{array}{ccccc}
0 & 0 & 0 & 0 & 0 \\
0 & 0 & 0 & 0 & x \\
0 & 0 & 0 & y & 0 \\
0 & 0 & y & a_{4,4} & a_{4,5} \\
0 & x & 0 & a_{4,5} & a_{5,5} 
\end{array}\right),$$
le lieu d{\'e}fini par le d{\'e}terminant est bien singulier ce qui montre que
l'on a un plan tangent {\`a} $\K^\vee$. On v{\'e}rifie ais{\'e}ment
que les fibres sont
bien form{\'e}es par les orbites sous l'involution (resp. les
involutions).\hfill$\Box$

\vs 0.2 cm

Choisir une quadrique $Q$ contenant $Y$ revient {\`a} choisir un point $x$
de $\p_3^\vee$. Lorsque $Y$ est lisse, 
la quadrique $Q$ d{\'e}finit un point
de  $U(A)$ si et seulement si elle ne contient aucune des 16 droites
$(p_iq_i)_{i\in[1,16]}$ avec $2(p_i+q_i)=\L-K_Y$.

La quadrique $Q$ peut {\^e}tre singuli{\`e}re. Nous expliquons ces cas. 
Identifions $A={\rm Jac}(Y)$ avec $\pic^6(Y)$ en choisissant le
diviseur $3K_Y$ pour origine. Ceci fixe 
un plan $H'\subset\p_3^\vee$ tangent {\`a} $\K^\vee$ selon une
conique. Soit $x\in U(A)\subset\p_3^\vee$ un relev{\'e} de
$\overline{x}\in U(A)/H$ d{\'e}finissant une
courbe de genre 5 et soit $\L\in\pic^6(Y)$. Alors, si l'image de $\L$
dans $\K$ n'est pas contenue dans $x^\vee$, la quadrique d{\'e}finie par
$Q$ est lisse, sinon elle est singuli{\`e}re. La quadrique est donc
singuli{\`e}re exactement si $\L\in C\subset A$. Nous verrons que ces cas
donnent des plongements particuliers de $C$. Remarquons que la courbe
$Y$ ne passe jamais par le sommet de $Q$ car sinon $Q$ d{\'e}finirait un
point de $\K^\vee$ ({\`a} l'intersection du plan tangent et de $\K^\vee$,
cf. proposition 3.1) et donc pas un point de $U(A)$.


Dans le second cas, 
notons $(x_i)_{i\in[1,4]}$ et $(y_j)_{j\in[1,4]}$ les points fixes des
involutions $i_1$ et $i_2$.
La quadrique $Q$ d{\'e}finit un point de $U(A)$ si et
seulement si $Q$ est de rang 4 et ne contient aucune droite $(x_iy_j)$. 

\vs 0.2 cm

Nous 
d{\'e}finissons la courbe $\widetilde{X}$ comme
l'image de $C$ dans $\K$ et la courbe $X$ comme l'image r{\'e}ciproque
de $\widetilde{X}$
dans $\widetilde{\K}$. Supposons maintenant que $Q$ d{\'e}finit un point
de $U(A)$, alors 

\vs 0.2 cm

\vs 0.2 cm

\thg{P}{roposition}{3.7} \textit{La courbe $C$ est lisse de genre 5 et
  c'est un rev{\^e}tement double non ramifi{\'e} de la courbe $X$ qui est
  lisse de genre 3.}

%

\subsection{Des plongements particuliers de la courbe de genre 5}

Nous avons vu que la donn{\'e}e d'une surface ab{\'e}lienne $A$ et d'un
point $\overline{x}\in U(A)/H$
correspond {\`a} la donn{\'e}e d'une courbe $C$ de genre 5
munie d'une involution sans point fixe.


Nous montrons maintenant que la donn{\'e}e du diviseur $\L$ de degr{\'e} 6 sur
$Y$ pr{\'e}c{\'e}dent d{\'e}finit un morphisme particulier de $C$ vers un espace
projectif de dimension 3. En effet, dans $\p(H^0\L)$ le point $x\in
U(A)$ correspond {\`a} une quadrique $Q$ (lisse ou de rang 4) contenant
$Y$. 
%
%
Si $Q$ est lisse, la vari{\'e}t{\'e} de ses droites forme un espace projectif de dimension 3 que nous notons $\p(V)$.

Si $Q$ est de rang 4, les droites de $Q$ forment une fibration en $\p_2$ au dessus de $\pu$. Il existe deux
  contractions naturelles de cette fibration vers un espace projectif
  de dimension 3 que nous notons $\p(V_k)$ pour $k\in\{1,2\}$. Elles
  contractent un $\pu\times\pu$ vers une droite $L_k$ selon l'une ou
  l'autre des projections. La quadrique singuli{\`e}re peut alors {\^e}tre vue
  comme la vari{\'e}t{\'e} des droites de $\p(V_k)$ qui rencontrent
  $L_k$. 

Lorsque $Y$ est lisse (cf. le cas 3.2.2, $(a,b)=(2,4)$), on note
  $\p(V)$ l'espace projectif $\p(V_k)$ tel que les $(\a)$-plans
  rencontrent $Y$ en 2 points.

Lorsque $Y$ est singuli{\`e}re, les
  plans $\Pi_1$ et $\Pi_2$ sont des $(\a)$-plans pour un des $\p(V_k)$
  et des $(\b)$-plans pour l'autre. Notons $\p(V)$ l'espace projectif
  $\p(V_k)$ tel que les plans $\Pi_i$ soient des $(\a)$-plans.


\vs 0.2 cm

Dans tous les cas on a donc un morphisme de $C$ vers un espace projectif
$\p(V)$ de dimension~3. 

\vs 0.2 cm

\thg{P}{roposition}{3.8} (\i) \textit{Lorsque $Y$ et $Q$ sont lisses, la
  courbe $C$ est plong{\'e}e dans $\p(V)$, elle est de degr{\'e} 8 et a une
  infinit{\'e} de quadris{\'e}cantes.}

(\i\i) \textit{Lorsque $Y$ est lisse mais $Q$ de
  rang 4, la courbe $C$ est plong{\'e}e dans $\p(V)$ et est de degr{\'e} 7. Il
  existe une droite $L$ telle que la courbe $C\cup L$ a une
  infinit{\'e} de quadris{\'e}cantes.}

(\i\i\i) \textit{Lorsque $Y$ est r{\'e}union de deux courbes elliptiques,
  l'image $\overline{C}$ de la courbe $C$ dans $\p(V)$ est de degr{\'e} 8
  et de genre arithm{\'e}tique 7. Elle est li{\'e}e {\`a} une droite $L$ par une
  intersection de deux c{\^o}nes cubiques. Elle a deux points doubles
  situ{\'e}s aux sommets des c{\^o}nes.}

\vs 0.2 cm

\dm
(\i) 
La courbe $C$ est d{\'e}finie comme
les points $\xi\in A$ tels que la bis{\'e}cante $L_\xi$ {\`a} $Y$
d{\'e}finie par $\xi$ est contenue dans la quadrique $Q$. Nous
d{\'e}finissons donc le morphisme $C\to\p(V)$ par 
$\xi\mapsto L_\xi$ 
Ce morphisme est {\`a} priori d{\'e}fini lorsque
$\xi\neq K_Y$ mais nous pouvons le prolonger, si besoin est, en ce point
car $C$ est une courbe lisse (ce cas appara{\^\i}t lorsque $Q$ est singuli{\`e}re).

Soit $p\in Y$, ce point d{\'e}finit une droite $\ell_p$ de $\p(V)$. Cette droite
est quadris{\'e}cante {\`a} l'image de $C$ dans $\p(V)$. En effet, on cherche
les points $q\in Y$ tels que la droite $(pq)$ soit contenue dans
$Q$. Il faut donc que $q$ soit dans l'orthogonal de $p$ pour la
quadrique. Cet hyperplan rencontre $C$ en 6 points et est tangent {\`a}
$C$ en $p$. Il reste donc 4 autres points qui donnent 4 points de
l'image de $C$ dans $\p(V)$ qui sont situ{\'e}s sur $\ell_p$.

La courbe $Y$ d{\'e}finit dans $\p(V)$ une surface r{\'e}gl{\'e}e (l'ensemble des
droites $\ell_p$ pour $p\in Y$). L'image de $C$ dans $\p(V)$ est
exactement le lieu singulier de cette surface. En effet, soit
$\xi\in C$, on a $\xi=p+q$ et l'image de $\xi$ dans $\p(V)$ est {\`a}
l'intersection des droites $\ell_p$ et $\ell_q$. R{\'e}ciproquement si un
point $z$ de $\p(V)$ est dans le lieu singulier de la surface, alors il
est sur deux droites $\ell_p$ et $\ell_q$, la droite de $\p_4$ que $z$
d{\'e}finit est alors contenue dans $Q$ et coupe $C$ en $p$ et $q$. Le
point $z$ est donc dans l'image de $C$.

Nous savons alors que la courbe image de $C$ dans $\p(V)$ est de degr{\'e}
8 et de genre arithm{\'e}tique 5 (cf. \cite{KL}). On a donc bien un
plongement qui est de degr{\'e} 8. La vari{\'e}t{\'e} des quadris{\'e}cantes de $C$
dans ce plongement est la courbe $Y$.

(\i\i) Notons $\overline{C}$ l'image de $C$. La quadrique $Q$ est
l'ensemble des droites rencontrant une droite $L$. Soit $S$ la surface
r{\'e}gl{\'e}e d{\'e}finie par $Y$ et $C'$ son lieu singulier. Un point $y$ de $L$
d{\'e}finit un$(\a)$-plan contenu dans $Q$. Il rencontre $Y$ en deux
points $y_1$ et $y_2$ en involution. Ainsi le point $y$ est {\`a}
l'intersection des droites $\ell_{y_1}$ et $\ell_{y_2}$. La droite $L$
est donc contenue dans $C'$. Par ailleurs, le m{\^e}me raisonnement que
pr{\'e}cedement montre que $\overline{C}\subset C'$ et que les points de $C'$ sont
donn{\'e}s par les bis{\'e}cantes de $Y$ contenues dans $Q$ (la diff{\'e}rence par
rapport {\`a} la situation pr{\'e}c{\'e}dente est le fait que le point $\xi=K_Y$
d{\'e}finit un $\pu$ --- qui donnera $L$ --- de bis{\'e}cantes
suppl{\'e}mentaires). On a donc $C'=\overline{C}\cup L$. Le genre de
$\overline{C}$ est sup{\'e}rieur ou {\'e}gal {\`a} celui de $C$ et par les
r{\'e}sultats de \cite{KL}, le genre de $C'$ est 5. On a deux cas :
$\overline{C}$ est de genre 6 et ne rencontre pas $L$ ou
$\overline{C}$ est de genre 5 et rencontre $L$ en un point.  La courbe
$\overline{C}$ rencontre $L$ : en effet, le point $K_Y$ est limite de
points $\xi\in C$ donc la limite des droites $L_\xi$ est une droite
d{\'e}finie par $K_Y$ donc une droite qui rencontre $L$. La courbe
$\overline{C}$ est donc de genre 5, c'est la courbe $C$. La courbe $C'$
est isomorphe {\`a} l'image r{\'e}ciproque de $C\subset J(Y)$ dans $\pic^2(Y)$
(qui est l'{\'e}clatement du point $K_Y$). Remarquons que dans cette
situation, les prongements sont obtenus {\`a} partir du plongement
canonique $K_C$ par projection {\`a} partir d'un point de $C$. Tous les
points de $C$ donnent un tel plongement. 

(\i\i\i) 
En d{\'e}crivant la quadrique comme la vari{\'e}t{\'e} des droites de
$\p(V)$ qui rencontrent un droite $L$, on voit que la courbe
$\overline{C}$ est sur l'intersection des deux c{\^o}nes cubiques
d{\'e}finis par les courbes $E_i$. Ces deux c{\^o}nes ont une droite en commun
car les courbes $E_i$ se coupent en un point. La courbe r{\'e}siduelle
doit alors {\^e}tre de degr{\'e} 8 et de genre arithm{\'e}tique 7. Elle a deux
points doubles. Chaque point double correspond aux points de $C$ de la
forme $(x,P)$ et $(P,y)$ ($P$ est le point d'intersection de $E_1$ et
$E_2$). Les deux points doubles sont donc aux sommets des
c{\^o}nes\hfill$\Box$

\vs 0.2 cm

\thf{Remarque}{3.9} Lorsque $Y$ est lisse et $Q$ lisse (resp. de rang
4), la courbe $C$ (resp $C'$) plong{\'e}e dans $\p(V)$ n'est pas 
contenue dans une surface cubique. En effet, supposons que c'est le
cas et soit $S$ une telle surface. Alors $S$ rencontre toutes les
quadris{\'e}cantes de $C$ (resp. $C'$) en quatre points donc $S$ contient
toutes les quadris{\'e}cantes de $C$ (resp. $C'$). La surface $S$ contient
donc la surface r{\'e}gl{\'e}e d{\'e}finie par $Y$, cette derni{\`e}re est de degr{\'e} 6,
c'est absurde.

\vs 0.2 cm

Nous montrons maintenant une r{\'e}ciproque {\`a} cette proposition dans
chacun des deux cas :

\vs 0.2 cm

\thg{P}{roposition}{3.10} \textit{Soit $\overline{C}$ une courbe de
  $\p(V)$ de degr{\'e} 8 et de genre g{\'e}om{\'e}trique 5 qui est li{\'e}e {\`a} une
  droite par une intersection de deux c{\^o}nes cubiques et qui a deux
  points doubles aux sommets des c{\^o}nes, alors le mod{\`e}le non
  singulier $C$ de $\overline{C}$ est dans $\M_5^{i,d}$.}

\vs 0.2 cm

\dm
Soit $L$ la droite comune aux deux c{\^o}nes cubiques. Soient $E_1$ et
$E_2$ les deux courbes elliptiques de la grassmannienne qui
d{\'e}finissent ces c{\^o}nes. La droite $L$ correspond au point
d'intersection des $E_i$.

Consid{\'e}rons l'incidence point/plan suivante :
$$I=\left\{(p,h)\in\p(V)\times\p(V)^\vee\ /\ p\in h\ {\rm et}\
  L\subset h\right\}.$$
La vari{\'e}t{\'e} 
%
$I$ param{\`e}tre les droites de la quadrique
singuli{\`e}re 
$$Q=\{l\in\G(2,\p(V))\ /\ l\cap L\neq\emptyset\}.$$ 
Notons $A=E_1\times E_2$, pour tout $\xi\in A^\times$, on peut d{\'e}finir un
{\'e}l{\'e}ment $L_\xi$ de $I$ comme {\'e}tant la droite passant par les deux
points de $\xi$. On a vu que la courbe
$$C=\overline{\{\xi\in A^\times\ /\ L_\xi\subset Q\}}$$ 
est alors lisse de genre 5. Elle s'envoie sur
$\overline{C}$ et on peut utiliser la construction pr{\'e}c{\'e}dente.\hfill$\Box$

\vs 0.2 cm

\thg{P}{roposition}{3.11} \textit{Soit $C$ une courbe projective de
  genre~5 et de degr{\'e}~8 de $\p_3$. Supposons que $C$ est lisse,
  irr{\'e}ductible et non hyperelliptique ou r{\'e}union d'une courbe lisse $C'$
de degr{\'e} 7 et d'une droite $L$ la rencontrant en un
  point. On suppose que $C$ admet une infinit{\'e} de droites
  quadris{\'e}cantes, alors $C\in\M_5^{i,i}$.}

\textit{La courbe $Y$ des quadris{\'e}cantes de $C$ est lisse de genre 2 et de
  degr{\'e} 6 (dans la grassmannienne). La courbe $C$ peut {\^e}tre retrouv{\'e}e
  {\`a} partir de la courbe $Y$ gr{\^a}ce {\`a} la construction pr{\'e}c{\'e}dente.}

\vs 0.2 cm

\dm
Remarquons que comme $C'$ est lisse de degr{\'e} 7 et de genre 5, elle est alors
non trigonale (cf \cite{Hartsh} exemple 6.4.2).

Nous commen\c cons par remarquer que $C$ ne peut {\^e}tre contenue dans
une surface cubique. En effet, soit $S$ une telle surface, la surface
$S$ rencontre toutes les quadris{\'e}cantes de $C$ en au moins 4
points. Elle contient donc toute les quadris{\'e}cantes de $C$. La surface
$S$ contient ainsi la surface r{\'e}gl{\'e}e $S'$ des quadris{\'e}cantes de
$C$. 
Le surface $S'$ ne peut {\^e}tre un plan. Si la surface
$S'$ {\'e}tait une quadrique lisse, alors $C$ serait lisse et
hyperelliptique.

Si $S'$ est un c{\^o}ne, alors $C$ passe par le sommet si et seulement si
elle est singuli{\`e}re. Mais alors $C'$ est trigonale (projection par le
sommet de degr{\'e} 3 vers $\pu$). Elle est donc lisse et ne passe pas par
le sommet. La
projection {\`a} partir de celui-ci donne un morphisme de degr{\'e} 4 de $C$
vers une conique plane. Les quadris{\'e}cantes coupent donc $C$ en
exactement 4 points. On rel{\`e}ve alors $C$ dans le mod{\`e}le non singulier
de $S'$ (l'{\'e}clatement du c{\^o}ne en son sommet qui est une surface
rationnelle $F_2$). On {\'e}crit $[C]=af+bh$ o{\`u} $f$ est la classe d'une
fibre et $h$ est la classe d'une section, voir par exemple
\cite{Beauville}. On a les intersections $f^2=0$, $f\cdot
h=1$ et $h^2=2$. Comme $C$ ne passe pas par le sommet du c{\^o}ne, sa
transform{\'e}e stricte ne rencontre pas le diviseur exceptionel $E$ dont
la classe est $h-2f$. On a donc $a=[C]\cdot[E]=0$. Par ailleurs $C$
rencontre les fibres (les quadris{\'e}cantes en exactement 4 points donc
$b=[C]\cdot f=4$. On a donc $[C]=4h$ ce qui donne par la formule
d'adjonction
$$2g(C)-2=[C]\cdot([C]+K)=4h(4h-2h)=16.$$
Ceci donne $g(C)=9$, c'est absurde.

Nous pouvons donc supposer que $S$ est la surface r{\'e}gl{\'e}e des
quadris{\'e}cantes. Le th{\'e}or{\`e}me de Segre (cf. \cite{GP2}
p.412) nous dit alors que si $\sigma$ est le genre de la courbe de
base de $S$ et si $C$ (resp. $C'$) rencontre les r{\`e}gles de $S$ avec
multiplicit{\'e} $k$, alors on a
$$2\ \!{\rm g}(C)-2=(k-1)(2\ \!{\rm deg}(C)-k\ \!{\rm deg}(S))+k\
\!(2\ \!\sigma-2)$$ 
et la m{\^e}me formule avec $C'$. Le genre $\sigma$ est $0$ ou $1$ ce qui donne 
$$3k^2-17k+24=0\ \ \ \ {\rm ou }\ \ \ \ 3k^2-19k+24=0\ \ {\rm pour}\
C\ {\rm lisse}$$
$$\ {\rm et}\ 3k^2-15k+22=0\ \ \ \ {\rm ou }\ \ \ \
3k^2-17k+22=0\ {\rm pour}\ C'.$$
La seule solution enti{\`e}re pour $C$ lisse est pour $\sigma=0$, on a
alors $k=3$ ce qui est impossible car les r{\`e}gles sont des
quadris{\'e}cantes {\`a} $C$ donc $k\geq 4$. Pour $C'$ la seule solution
enti{\`e}re est pour $\sigma=1$ ce qui donne $K=2$ ce qui est encore
impossible, les r{\`e}gles sont des trid{\'e}cantes {\`a} $C'$.

\vs 0.2 cm

Nous montrons maintenant le

\vs 0.2 cm

\thg{L}{emme}{3.12} \textit{La courbe n'a pas de quintis{\'e}cante et il
  existe une surface de degr{\'e} inf{\'e}rieur ou {\'e}gal {\`a} sept singuli{\`e}re le
  long de $C$.}

\vs 0.2 cm

\dm
On peut mettre une structure de vari{\'e}t{\'e} sur la famille des
quadris{\'e}cantes (cf. \cite{GP} ou \cite{ACGH}). L'hypoth{\`e}se signifie
qu'il existe une courbe de la grassmanienne des droites de
$\mathbb{P}_3$ correspondant {\`a} des quadris{\'e}cantes de $C$.
La courbe $C$ n'a pas de 5-s{\'e}cante. En effet, si $D$
est une 5-s{\'e}cante, on a une suite exacte
$$0 \to \oo_C(-Z) \to \oo_{C\cup D} \to \oo_D \to 0$$
o{\`u} $Z = C\cap D$. Il en r{\'e}sulte la suite exacte longue
$$0 \to H^0(\oo_C(3h - Z)) 
\to H^0(\oo_{C\cup D}(3)) 
\to H^0(\oo_D(3))\to 0$$
o{\`u} $h$ est la classe d'une section plane de $C$ et donc que
$H^0(\oo_{C\cup D}(3)) = 19$.
On en d{\'e}duit que $C\cup D$ est sur une surface cubique ce qui est
absurde. 

Il existe une surface de degr{\'e} inf{\'e}rieur ou {\'e}gal {\`a} 7 singuli{\`e}re le long de
$C$. En effet, comme $C$ n'est pas sur une 
cubique, on a
$h^1\I_C(3)=0$. De plus 
$h^2\I_C(2)=h^1\oo_C(2)=0$ et
$h^3\I_C(1)=0$. Le faisceau $\I_C$ est donc 3-r{\'e}gulier ce qui impose
que $\I_C(4)$ est engendr{\'e} par ses sections (cf \cite{mumcurves}). On
a donc 
une surjection $\oo_C(3)\to(\I_C/\I_C^2)(7)$ et comme
$h^1\oo_C(3)=0$, on en d{\'e}duit que $h^1((\I_C/\I_C^2)(7))=0$ et donc
$h^0((\I_C/\I_C^2)(7))=\chi((\I_C/\I_C^2)(7))=64$. On a donc
$h^0(I^2_C(7))
\ge4$.
Il existe donc une surface de degr{\'e} inf{\'e}rieur ou {\'e}gal {\`a} 7
singuli{\`e}re le long de $C$. Cette surface contient toutes les
quadris{\'e}cantes de~$C$.\hfill$\Box$

%

\vs 0.2 cm

\thg{L}{emme}{3.13} \textit{Si deux quadris{\'e}cantes $l$ et $l'$ de $C$
  se coupent en un point $p$ alors $p$ appartient {\`a} $C$. Les droites
  $l$ et $l'$ et la tangente {\`a} $C$ en $p$ ne sont pas coplanaires.}

\vs 0.2 cm

\dm 
Si $p\not\in C$, on a la suite exacte 
$$0\to \oo_C (-(C\cap (l \cup l')) \to \oo_{C\cup l\cup l'}
\to \oo_{l\cup l'} \to 0.$$

Apr{\`e}s tensorisation par $\oo_{\p_3}(3)$, on obtient $h^0(\oo_C(3h -
p_1 - \ldots - p_4-q_1 \ldots - q_4) = 12$ et $h^0(\oo_{l\cup l'}(3))
= 7$. Il en r{\'e}sulte que $h^0(\oo_{C\cup l\cup l'}(3)) = 19$
donc $C\cup l\cup l'$ est sur une surface cubique ce qui est absurde. 
On a ici not{\'e} $l\cdot C = p_1 + p_2 + p_3 + p_4$ et $l'\cdot C
= q_1 + q_2 + q_3 + q_4$. 

Si le plan engendr{\'e} par $l$ et $l'$ est tangent {\`a} $C$ en $p$,
toute surface contenant $l$ et $l'$ est tangente {\`a} $C$ en $p$
et on conclut comme pr{\'e}c{\'e}demment avec $l\cdot C = p+p_1+p_2+p_3$ et
$l'\cdot C = p+q_1+q_2+q_3$ et en consid{\'e}rant $h^0\big(\oo_C (3h-2p-
p_1-p_2-p_3-q_1-q_2-q_3)\big)$. On peut, de la m{\^e}me mani{\`e}re, montrer
que par un point de $C$ il ne passe pas trois quadris{\'e}cantes.\hfill$\Box$

\vs 0.2 cm

Soit $\bar S$ la surface r{\'e}gl{\'e}e engendr{\'e}e par la courbe $Y$ des
quadris{\'e}cantes.

\vs 0.2 cm

\thg{L}{emme}{3.14} \textit{La courbe $C$ est le lieu singulier de la
  surface $\bar S$ qui est de degr{\'e} 6. La courbe $Y$ est lisse et de
  genre 2.}

\vs 0.2 cm

\dm
Le degr{\'e} de $\bar S$  est inf{\'e}rieur ou {\'e}gal {\`a} 7
d'apr{\`e}s le lemme 3.12 et est sup{\'e}rieur ou {\'e}gal {\`a} quatre. Soit $Y$ la
courbe des quadris{\'e}cantes. Elle est irr{\'e}ductible sinon l'une de ses
composantes donne une surface de degr{\'e} inf{\'e}rieur {\`a} 3 contenant
$C$. Remarquons que si $Y$ est plane alors $\bar S$ est un c{\^o}ne et soit
$x$ son sommet. Si $C$ est lisse et $x\not\in C$, alors la projection
{\`a} partir de $x$ est un morphisme de degr{\'e} 4 (car les quadris{\'e}cantes
passent par $x$), son image est donc de degr{\'e} 2 et $C$ est sur un
c{\^o}ne, c'est impossible. Si $x\in C$ et $C$ lisse, alors le morphisme
de projection doit {\^e}tre degr{\'e} au moins 3 et diviser 7, il est donc de
degr{\'e} 7, son image est alors de degr{\'e} 1 et $C$ serait plane,
absurde. Si $C$ n'est pas lisse et $x\not\in C'$ alors la projection
de $C'$ {\`a} partir de $x$ est de degr{\'e} au moins 3 (car les droites sont
au moins tris{\'e}cantes), elle est donc de degr{\'e} 7 ce qui impose que $C'$
est plane, c'est absurde. Si $x\in C'$, la projection est encore au
moins de degr{\'e} 3. Si elle est de degr{\'e} 3 alors on a un ${\bf g}_3^1$
sur $C'$, c'est absurde. Elle doit donc {\^e}tre de degr{\'e} 6 et $C'$ plane,
c'est encore impossible. Ainsi $Y$ n'est pas plane.

Le lieu singulier de $\bar S$ est contenu dans $C$ d'apr{\`e}s
le lemme pr{\'e}cedent, il est donc de degr{\'e} 0, 1, 7 ou 8. Soit $n$ le
degr{\'e} de $Y$ et $p$ son genre g{\'e}om{\'e}trique. Les r{\'e}sultats de \cite{KL}
nous disent que le lieu singulier de $\bar S$ est alors de degr{\'e}
$(n-1)(n-2)/2 - p$ qui vaut 
$3-p$, $6-p$, $10-p$ ou $15-p$ selon les valeurs de
$n\in[4,7]$. De plus ce degr{\'e} doit appartenir {\`a} $\{0,1,7,8\}$. Mais
comme $Y$ n'est pas contenue dans un plan et est irr{\'e}ductible la
formule de Castelnuovo (cf. \cite{ACGH} p.116) nous donne une borne
sur $p$ qui est 1, 2, 4 ou 6 selon les valeurs de $n$. Ainsi
le degr{\'e} du lieu singulier est contenu dans $[2,3]$, $[4,6]$,
$[6,10]$ ou $[9,15]$ et dans $\{0,1,7,8\}$. Les seules valeurs
possibles sont 7 ou 8 avec $n=6$ et $p=3$ ou 2.

Soit $\widetilde{Y}$ le mod{\`e}le lisse de $Y$. Le lieu singulier de
$\bar S$ peut {\^e}tre vu comme une courbe dans ${\rm
  Div_2}(\widetilde{Y})$ dont le genre arithm{\'e}tique est $P$ avec $2P-2 =
(n-5)(n+2p-2)$ (cf. \cite{A}). Ceci nous donne $P=6$ si $p=3$ ou $P=5$ si
$p=2$. Le premier cas est impossible. Ainsi $Y$ est de degr{\'e} 6 et de
genre g{\'e}om{\'e}trique 2. On peut calculer (cf. \cite{KL}) que la surface n'a
pas de point triple.

Il reste {\`a} montrer que $Y$ est lisse. 
Si $Y$ n'est pas lisse, alors $p_a(Y)\geq 3$ et 
$Y$ est contenue dans l'intersection de la grassmannienne
avec un espace projectif de dimension 3. Si cette intersection est
lisse c'est une quadrique form{\'e}e des droites rencontrant deux droites
$L_1$ et $L_2$. Les deux droites sont alors dans le lieu singulier,
c'est impossible. Si l'intersection est singuli{\`e}re, on peut supposer
que c'est un c{\^o}ne (car $Y$ n'est pas plane). Mais alors un calcul dans
la surface $F_2$ mod{\`e}le non singulier de ce c{\^o}ne montre que le genre
g{\'e}om{\'e}trique de $Y$ doit {\^e}tre 3 ou 4. C'est impossible donc $Y$ est
lisse.\hfill$\Box$

\vs 0.2 cm

\textit{Preuve de la proposition 3.11 :}
Notons $\G$ la grasmannienne des droites de $\p_3$, on voit $\G$ comme
une quadrique de $\p_5$ par le plongement de Pl{\"u}cker. Notons $E$ la
restriction du quotient tautologique de $\G$ {\`a} $Y$. Le mod{\`e}le non
singulier de $\bar S$ est $\p_Y(E)$.

Soit $\L=\Lambda^2E=\oo_\G(1)\vert_C$, c'est un fibr{\'e} de degr{\'e} 6. L'espace
$\p_4=\p(H^0\L)$ est un hyperplan de $\p_5$. La courbe $Y$ est
contenue 
dans la quadrique $Q$ d{\'e}coup{\'e}e par $\G$
dans $\p_4$. 

La courbe $C$ est isomorphe {\`a} la courbe suivante de $Y_2$ le carr{\'e}
sym{\'e}trique de $Y$ :
$$\tilde{C}=\{(x,y)\in Y_2\ /\ (xy)\subset Q\}.$$
En effet, on d{\'e}finit un morphisme $\tilde{C}\to C$ par $(x,y)\mapsto\ell_x\cap
\ell_y$ o{\`u} $\ell_x$ (resp. $l_y$) d{\'e}signe la droite de $\p_3$ correspondant
{\`a} $x$ (resp. $y$). Ce morphisme est bien d{\'e}fini car si on a
$(xy)\subset Q$, alors les droites $\ell_x$ et $\ell_y$ se rencontrent
et on a (d'apr{\`e}s le lemme 3.13) un point de $C$.
On d{\'e}finit une r{\'e}ciproque $C\to\tilde{C}$ par $z\mapsto(x,y)$ o{\`u} $x$
et $y$ sont les deux quadris{\'e}cantes passant par $z$ (il y en a
exactement deux car la surface n'a pas de point triple). 


Nous sommes dans la situation de la proposition 3.8. En particulier on
a une involution $i$ sur $\tilde{C}$ (vue dans $\pic_2(Y)$) donn{\'e}e par
$\xi\mapsto\L-K_Y-\xi$ (cf. proposition 3.1). 
Elle est
donn{\'e}e sur $C$ de la mani{\`e}re suivante : soit $z\in C$, il existe deux
quadris{\'e}cantes $L_1$ et $L_2$ passant par $z$. Soit alors $H$ le plan
engendr{\'e} par ces quadis{\'e}cantes, le plan $H$ rencontre $C$ en 8 points
dont 7 sont sur $L_1\cup L_2$. Le huiti{\`e}me point est
$i(z)$. L'involution est sans point fixe. En effet, si $z=i(z)$, ceci
signifie que $H$ est tangent {\`a} $C$ en $z$ et $H$ contient $L_1$ et
$L_2$. C'est exclu par le lemme 3.13.\hfill$\Box$

\vs 0.2 cm

\thf{Remarque}{3.15} Dans \cite{M}, 
S. Mukai donne une version du lemme~3.12. 

\vs 0.2 cm

Ces r{\'e}sultats nous permettent de donner une nouvelle caract{\'e}risation
g{\'e}om{\'e}trique des cour\-bes de genre 5 rev{\^e}tement double d'une courbe de
genre 3 :

\vs 0.2 cm

\thg{T}{h{\'e}or{\`e}me}{3.16} \textit{Soit $C$ une courbe lisse de genre
  5.}

(\i) \textit{{\rm Premi{\`e}re caract{\'e}risarion :} On a $C\in\M_5^{i,i}$ si
  et seulement s'il existe un plongement ${\cal M}$  de degr{\'e} 8 de $C$
  dans $\p_3$ pour lequel la courbe $C$ a une infinit{\'e} de
  quadris{\'e}cantes.}

\textit{{\rm Seconde caract{\'e}risarion :} On a $C\in\M_5^{i,i}$ si
  et seulement s'il existe un plongement ${\cal M}'$  de degr{\'e} 7 de $C$
  dans $\p_3$ et une droite $L$ rencontrant $C$ en un point pour lesquels la courbe $C\cup L$ a une infinit{\'e} de quadris{\'e}cantes.}


\vs 0.1 cm

\textit{Dans cette situation, la courbe $Y$ des quadris{\'e}cantes est la m{\^e}me
  quelque soit le plongement, elle est lisse de genre 2 et telle que
  $J(Y)={\rm Prym}(C)$. Notons
  $\mathfrak{J}_0$ les diviseurs ${\cal M}$ de degr{\'e} 8 du premier
  type, $\mathfrak{J}_1$ les diviseurs ${\cal M}'$ de degr{\'e} 7 du
  second type et $\mathfrak{J}$ la r{\'e}union de ces ensembles.  Il y a
  un morphisme $\mathfrak{J}\to J(Y)$ qui est un fibr{\'e} principal
  homog{\`e}ne de groupe $H$.}

\vs 0.1 cm

(\i\i) \textit{On a $C\in\M_5^{i,d}$ si et seulement s'il existe un
  morphisme de $C$ dans $\p_3$ tel que son image $\overline{C}$ est
  li{\'e}e {\`a} une droite $L$ par une intersection de deux c{\^o}nes cubiques,
  de degr{\'e} 8 et a deux points doubles aux sommets des c{\^o}nes.}

\vs 0.1 cm

\textit{Dans cette situation, notons $\mathfrak{J}$ les diviseurs de
  degr{\'e} 8 de $C$ qui d{\'e}finissent de tels morphismes. On a un morphisme
  $\mathfrak{J}\to E_1\times E_2$ qui est un fibr{\'e} principal
  homog{\`e}ne de groupe $H$ (les $E_i$ sont les courbes elliptiques
  d{\'e}finissant les c{\^o}nes cubiques).}

\vs 0.2 cm

\dm
Les caract{\'e}risations de $\M_5^{i,i}$ et $\M_5^{i,d}$ d{\'e}coulent des
propositions 3.8, 3.10 et 3.11. Nous d{\'e}crivons les ensembles
$\mathfrak{J}$ de plongements.

(\i) Pla\c cons nous dans le cas $C\in\M_5^{i,i}$ et notons $A=J(Y)$ la
vari{\'e}t{\'e} de Prym associ{\'e}e. la proposition 3.11 nous dit que si on a un
plongement ${\cal M}\in\mathfrak{J}$ de la courbe $C$, alors la courbe
$Y$ est munie d'un plongement $\L$ de degr{\'e} 6 (le plongement de
Pl{\"u}cker) et d'une quadrique $Q$ la contenant (la grassmannienne). La
courbe $C$ fixe un point $\overline{x}\in U(A)/H$ et la quadrique $Q$
d{\'e}finit un rel{\`e}vement $x\in U(A)$ de ce point. On a donc un morphisme 
$$\mathfrak{J}\to A$$
\vs -0.9 cm
$${\cal M}\mapsto\L.\ \ $$
La proposition 3.8 nous dit que ce morphisme est surjectif et que sa
fibre est donn{\'e}e par les quadriques $Q$ qui rel{\`e}vent $\overline{x}$,
c'est-{\`a}-dire par le groupe $H$ des {\'e}l{\'e}ments d'ordre 2 de $A$. 

Le morphisme $\mathfrak{J}\to A$ est le morphisme ${\rm
  Im}(1-i)\stackrel{1-i}{\fll}{\rm Im}(1-i)$ ou enore le morphisme
  $J(Y)\stackrel{2}{\fll}J(Y)$. Remarquons enfin que l'ensemble
  $\mathfrak{J}_1$ des diviseurs du second type est isomorphe {\`a} la
  courbe $C$ : ce sont les projections du plongement canonique de $C$
  par un point de $C$. $\mathfrak{J}_1$ est donc l'image du
  plongement de $C$ dans $J(Y)$.

Le m{\^e}me raisonnement donne le cas (\i\i).\hfill$\Box$

\vs 0.5 cm

\noi
\textsc{Jean D'Almeida}

\noi
Universit{\'e} des Sciences et Technologies de Lille

\noi
F-59665 Villeneuve d'Ascq Cedex, France.

\vs 0.5 cm

\noi
\textsc{Laurent Gruson}

\noi
Universit{\'e} de Versailles

\noi
45 avenue des {\'E}tats-Unis

\noi
F-78035 Versailles, France.

\vs 0.5 cm

\noi
\textsc{Nicolas Perrin}

\noi
Institut de Math{\'e}matiques de Jussieu

\noi
175 rue du Chevaleret

\noi
F-75013 Paris, France.

\end{document}